\documentclass[pdftex, noinfoline]{imsart}

\RequirePackage{amsthm,amsmath,amsfonts,amssymb}
\RequirePackage[numbers]{natbib}
\RequirePackage[colorlinks,citecolor=blue,urlcolor=blue]{hyperref}
\RequirePackage{graphicx}
\usepackage{bm}
\usepackage{bbm}
\usepackage{mathrsfs}
\usepackage{url}

\startlocaldefs
\theoremstyle{plain}

\theoremstyle{definition}

\theoremstyle{remark}

\endlocaldefs

\begin{document}
\begin{frontmatter}
\title{Bayesian inference in high-dimensional models}
	
%\title{A sample article title with some additional note\thanksref{t1}}
\runtitle{Bayesian inference in high-dimensional models}
%\thankstext{T1}{A sample additional note to the title.}

\begin{aug}
%%%%%%%%%%%%%%%%%%%%%%%%%%%%%%%%%%%%%%%%%%%%%%%
%% Only one address is permitted per author. %%
%% Only division, organization and e-mail is %%
%% included in the address.                  %%
%% Additional information can be included in %%
\author[A]{\fnms{Sayantan}~\snm{Banerjee}\ead[label=e1]{sayantanb@iimidr.ac.in}},
\author[B]{\fnms{Isma\"el}~\snm{Castillo}\ead[label=e2]{ismael.castillo@upmc.fr}}
\and
\author[C]{\fnms{Subhashis}~\snm{Ghosal}\ead[label=e3]{sghosal@ncsu.edu}}
%%%%%%%%%%%%%%%%%%%%%%%%%%%%%%%%%%%%%%%%%%%%%%
%% Addresses                                %%
%%%%%%%%%%%%%%%%%%%%%%%%%%%%%%%%%%%%%%%%%%%%%%
\address[A]{Operations Management \& Quantitative Techniques Area, IIM Indore, Rau-Pithampur Road, Indore, MP 453556, India\printead[presep={,\ }]{e1}}

\address[B]{Laboratoire de Probabilit\'es, Statistique et Mod\'elisation, Sorbonne Universit\'e, Paris, France\printead[presep={,\ }]{e2}}

\address[C]{Department of Statistics, North Carolina State University, Raleigh, NC, U.S.A.\printead[presep={,\ }]{e3}}

\runauthor{S.Banerjee et al.}
\end{aug}

\begin{abstract}
Models with dimensions exceeding the available sample size are now commonly used across various applications. A sensible inference is possible using a lower-dimensional structure. In regression problems with many predictors, the model is often assumed to be sparse, with only a few predictors active. Interdependence among a large number of variables is succinctly described by a graphical model, in which variables are represented as nodes on a graph, and an edge between two nodes indicates their conditional dependence given other variables. Many procedures for making inferences in the high-dimensional setting, typically using penalty functions to induce sparsity in the solution obtained by minimizing a loss function, were developed. Bayesian methods have been proposed more recently for such problems, where the prior accounts for the sparsity structure. These methods naturally quantify the uncertainty of the inference through the posterior distribution. Theoretical studies of Bayesian procedures in high dimensions have recently been conducted. Questions that arise are whether the posterior distribution contracts at the minimax optimal rate near the true parameter value, whether the correct lower-dimensional structure is discovered with high posterior probability, and whether a credible region has adequate frequentist coverage. In this paper, we review the properties of Bayesian and related methods for several high-dimensional models, including the many normal means problem, linear regression, generalized linear models, and Gaussian and non-Gaussian graphical models. Practical computational approaches are also discussed. 
\end{abstract}



\end{frontmatter}
%%%%%%%%%%%%%%%%%%%%%%%%%%%%%%%%%%%%%%%%%%%%%%
%% Please use \tableofcontents for articles %%
%% with 50 pages and more                   %%
%%%%%%%%%%%%%%%%%%%%%%%%%%%%%%%%%%%%%%%%%%%%%%
%\tableofcontents

	\section{Introduction}
	\label{sec:1}
	
	Advances in technology have produced massive datasets collected across all aspects of modern life. Huge datasets come from internet searches, mobile apps, social networking, cloud computing, wearable devices, as well as from more traditional sources such as bar-code scanning, satellite imaging, air traffic control, banking, finance, and genomics. Due to the complexity of such datasets, flexible models are needed involving many parameters, routinely exceeding the sample size. In such a situation, a meaningful inference is possible only if there is a hidden lower-dimensional structure involving far fewer parameters. This will occur in a linear or generalized linear regression model if the vector of regression coefficients is mostly zeros. In this situation, sensible inference is possible by forcing the estimated coefficient to be sparse via an automatic mechanism determined by the data.
    
Another critical problem is the study of the interrelationships among a large class of variables. It is often sensible to think that only a few pairs have some intrinsic relations between them, when the effects of other variables are eliminated; that is, most pairs of variables are conditionally independent given other variables. The underlying sparsity is conveniently described by a graph, where nodes represent variables, and an edge between two variables is present only if they are conditionally dependent given other variables. Hence, the resulting model is called a graphical model. When the variables are jointly Gaussian, the absence of an edge is equivalent to having a zero entry in the precision (inverse of the covariance) matrix. Thus, learning the structural relationship in such a Gaussian graphical model is possible by forcing the estimated precision matrix to have zero entries in most places. Other problems that effectively use a lower-dimensional structure include matrix completion problems (where many entries of a matrix are missing and it is assumed that the underlying true matrix has a sparse plus a low-rank structure), and stochastic block models (where the extent of interaction between two nodes is determined solely by their memberships in certain hidden blocks of nodes). 
	
	Numerous methods for estimating parameters in the high-dimensional setting have been proposed in the literature, most of which use penalization. The idea is to add a suitable penalty term to the loss function to be optimized, ensuring the resulting solution is sparse. The most familiar method is the least absolute deviation shrinkage and selection operator, abbreviated as LASSO (Tibshirani \cite{tib}): 
    \begin{align*}
    \hat \beta =\arg \min \sum_{i=1}^n (Y_i- \sum_{j=1}^p \beta_j X_{ij})^2+\lambda \sum_{j=1}^p |\beta_j| 
    \end{align*}
    for a linear regression model $Y_i=\sum_{j=1}^p \beta_j  X_{ij}+\varepsilon_i$, $i=1,\ldots,n$. The sharp corners of the contours of the $\ell_1$-penalty function $\sum_{j=1}^p |\beta_j|$ in this case force sparsity in the resulting estimate, where the extent of the sparsity depends on the tuning parameter $\lambda$ and the data. While there are many other important penalization procedures for this problem and many other related problems, we shall not explicitly refer to them except to introduce and compare with a Bayesian method, which is the primary object of interest in this review. {Excellent sources of information on the frequentist literature for high-dimensional models are B\"uhlmann and van de Geer \cite{buehlmannvdGeer2011}, Giraud \cite{giraudbook}, Wainwright \cite{wainwrightbook}}.  Bayesian methods for high-dimensional models have attracted significant recent interest. A hidden lower-dimensional structure, such as sparsity, can be easily incorporated into a prior distribution, for instance, by allowing a point mass (or some analog of it) at zero. In addition to providing a point estimate, a Bayesian method naturally assesses uncertainty in structure learning and provides credible sets with attached probabilities for quantifying uncertainty. Good Bayesian procedures should have desirable frequentist properties, such as a minimax-optimal rate of contraction, consistency in variable selection (or, more generally, structure learning), and asymptotic frequentist coverage of credible sets. Investigating these properties is a primary research objective, critical in high-dimensional settings, where identifying appropriate priors is more challenging. Continued research has shown that under suitable choices of prior distributions, Bayesian methods can asymptotically perform very well. 
	However, a Bayesian method is typically also very computationally intensive, especially if a traditional Markov chain Monte Carlo (MCMC) procedure is used, since MCMC iterations must explore many possible models. In recent years, effective computational strategies using continuous shrinkage priors, expectation-maximization (EM) algorithm, variational approaches, Hamiltonian Monte Carlo, and Laplace approximation have been proposed. Our review will address some aspects of computation in the Appendix. 
	
	We shall cover the following models of interest: Normal sequence model (Section~3); High-dimensional linear regression model (Section~4); High-dimensional nonparametric and other regression models (Subsection~4.3);  Gaussian graphical model (Subsection~5.2); High-dimensional classification (Subsection~5.3); Ising and other non-Gaussian graphical models (Subsection~5.4); Nonparanormal graphical models (Subsection~5.5); Signal on a large graph (Section~6); Matrix models such as structured sparsity and Stochastic block models (Section~7). We shall address the issues of posterior contraction, variable or feature selection, distributional approximation, and coverage of credible sets, as well as modifications to the Bayesian approach using fractional posterior and PAC-Bayesian methods. 
	
	{\it Notations}: A generic statistical model indexed by a parameter $\theta$ is written as $\mathcal{P}=\{P_\theta^{(n)},\ \theta\in\Theta\}$, with $\Theta$ a subset of a (typically, here, high-dimensional) metric space equipped with a distance $d$. 
    
    Let $\Pi$ stand for a prior distribution on $\Theta$, and suppose that one observes data $Y^{(n)} \sim P_\theta^{(n)}$.
	Let the true value of $\theta$ be denoted by $\theta_0$ so that we study convergence under the distribution $P_{\theta_0}$. We say that $\epsilon_n$ is a contraction rate (at $\theta_0$, with respect to a metric $d$ on $\Theta$) if $\mathrm{E}_{\theta_0}\Pi[ d(\theta,\theta_0)\le M_n\epsilon_n | X^{(n)}] \to 1$ as $n\to\infty$ for every sequence $M_n\to \infty$; see Ghosal and van der Vaart \cite{gvbook17}. 
    
    Let $\mathrm{N}(0,\Sigma)$ stand for the centered Gaussian distribution with covariance $\Sigma$ and $\mathrm{Lap}(\lambda)$ denote the Laplace distribution of parameter $\lambda>0$. Let  $I_d$ be the identity matrix of dimension $d\ge 1$, $\|\cdot\|_p$, $p\ge 1$, be the $\ell_p$-norm, $\|\cdot\|$ stand for the Euclidean i.e., $\|\cdot\|_2$-norm, and $\|\cdot\|_{\mathrm{F}}$ for the Frobenius (Euclidean) norm on matrices.

	\section{Priors for sparsity}
	\label{sec:sparse prior}
	
	\subsection{Spike-and-slab priors}
	
	As sparsity plays a vital role in high-dimensional modeling, sparsity-inducing priors are of particular interest. While an entry $\theta$ may likely be zero, the possibility of a non-zero value, even a large value, cannot be ruled out. This can be thought of as two regimes superimposed --- one corresponding to zero or minimal values of $\theta$, and the other to possibly large values --- and the prior should provide a probabilistic mechanism for choosing between the regimes. To address this, a {\em spike-and-slab prior} has been considered (Mitchell and Beauchamp \cite{mitchellbeauchamp}, George and McCulloch \cite{george93}, Ishwaran and Rao \cite{ishwaranrao}): 
	\begin{equation}
		\label{s&s density}
		\pi(\theta) =(1-w) \phi_0(\theta)+ w \phi_1(\theta), 
	\end{equation}
	where $\phi_0$ is a density highly concentrated at $0$, $\phi_1$ is a density (usually symmetric about $0$) allowing intermediate and large values of $\theta$, and $w$ is a small parameter, thus inducing sparsity in the mixture. For instance, $\phi_0$ may be the normal density with mean zero and a small variance $v_0$, and $\phi_1$ the normal density with mean zero and a relatively large variance $v_1$. The parameters in $\phi_0$ and $\phi_1$, as well as $w$, may be given further priors. Another choice, such that both $\phi_0$ and $\phi_1$ are Laplace densities, was proposed and called the  {\em spike-and-slab LASSO} (Ro\u ckov\`a and George \cite{rockovageorge18}). The primary motivation is to identify the posterior mode as a sparse vector, as in the LASSO. Generally, the spike part of the prior induces a shrinkage towards zero, which can be limited by using a heavier-tailed density $\phi_1$ for the slab, such as a t-density, or at least as heavy-tailed as the Laplace density. For the spike part, an extreme choice is the distribution degenerate at 0, which corresponds to exact sparsity, and the resulting prior will be referred to as the \emph{hard-spike-and-slab prior}, while the term \emph{soft-spike-and-slab prior} will be used if the spike has a density. Johnson and Rossell \cite{johnsonrossell10,johnsonrossell12} argued in favor of non-local priors,  which make the spike as separated as possible from the slab {around zero}, by choosing slab distributions that have very little mass close to $0$. George and McCulloch \cite{george93} developed the Stochastic Search Variable Selection (SSVS), a Gibbs sampling procedure, to compute the posterior distribution based on a soft-spike-and-slab prior, by augmenting indicator variables for slabs. Since the number of all possible models $2^p$ is enormous when the number of predictors $p$ is high, a direct evaluation of the posterior probabilities of each model is impossible, but the MCMC procedure is designed to visit models with relatively high posterior probabilities while passing the models with lower posterior probabilities. 
	
	\subsection{Continuous shrinkage priors}
	
	Computation using a spike-and-slab prior involves a latent indicator of the mixture component. Replacing the indicator by a continuous variable leads to the so-called {\em continuous shrinkage priors}, typically obtained as scale mixtures of normal. Alternative terms like global-local priors or one-component priors are also used. An early example is the Laplace prior (Park and Casella \cite{BayesianLasso}, Hans \cite{Hans}), which is an exponential scale-mixture of normals, often called the \emph{Bayesian LASSO} because the corresponding mode is interpreted as the LASSO. However, the obvious drawback of these priors is that there is not sufficient concentration near the value 0 for the entire posterior to concentrate near zero whenever a coefficient is 0, even though the posterior mode may be sparse or nearly sparse. The prior concentration near zero should be high while still maintaining a thick tail, by letting the scale parameter have a more spiked density at 0. A popular continuous shrinkage prior meeting these requirements is the {\em horseshoe} prior (Carvalho et al. \cite{carvalho}), which is a half-Cauchy scale mixture of normal: 
	\begin{equation} 
		\label{horseshoe}
		\theta | \lambda  \sim \mathrm{N}(0,\lambda^2),\qquad
		\lambda  \sim \mathrm{Cauchy}^+(0,\tau),
	\end{equation}
	where $\mathrm{Cauchy}^+(0,\tau)$ is the half-Cauchy distribution with scale $\tau$. The corresponding marginal density of $\theta$ has a pole at $0$ and Cauchy-like tails.
	A further prior may be put on $\tau$, such as another half-Cauchy prior leading to the Horseshoe+ prior (Bhadra et al. \cite{bhadra2017}).

	One may consider different priors on the scales $\lambda$ in \eqref{horseshoe} with a high concentration near 0. The scale parameter $\lambda$ is entry-specific and is called the local shrinkage prior. The scale parameter $\tau$ is common for all entries and is called the global shrinkage parameter. The local shrinkage parameter should have a heavy tail, while the global shrinkage parameter should have a high concentration at 0 (Polson and Scott \cite{Polson&Scott}), respectively controlling the tail and the sparsity. Various choices of mixing distributions lead to the introduction of many continuous shrinkage priors, such as normal-inverse-Gaussian prior (Caron and Doucet \cite{carondoucet}), normal-gamma prior (Griffin and Brown \cite{griffin2010}), and the generalized double Pareto priors (Armagan et al. \cite{armagan2013generalized}). Another possibility is the  Dirichlet-Laplace prior (Bhattacharya et al. \cite{bpd15}): 
	\begin{equation}  
		\label{Dirichlet-Laplace}
		\theta_i| \phi,\tau \sim \text{Lap}(\phi_i\tau),\qquad 
		\phi=(\phi_1,\ldots,\phi_p) \sim \text{Dir}(a,\ldots,a),
	\end{equation}
	where choosing $a\in(0,1)$ leads to a pole at $0$ for the marginal distribution of $\theta_i$'s enforcing (near-)sparsity, and  $\tau$ is a parameter which is given a gamma distribution.
	
	\section{Normal sequence model}
	\label{sec:2}
	
	The normal sequence model gives the simplest high-dimensional model: $Y_i=\theta_i+\varepsilon_i$, $i=1,\ldots,n$, where $\varepsilon_i$ are independent and identically distributed (i.i.d.) $\mathrm{N} (0,1)$, the parameter set $\Theta$ for $\theta=(\theta_1,\ldots,\theta_n)$ is $\mathbb{R}^n$. The posterior contraction rate will be obtained uniformly for parameters belonging to 
	a {\em nearly-black} class $\ell_0[s] = \left\{\theta\in\mathbb{R}^n: \#\{i: \theta_i\neq 0\}\le s\right\}$,  $0\le s=s(n)\le n$,  where $\#$ stands for the cardinality of a finite set.  {In this model, the minimax risk for the squared error loss is equivalent to} $2 s \log (n/s)$, see Donoho et al. \cite{dhjs92}. Let $\theta_0$ stand for the true value of the vector $\theta$ and $s_0$ the cardinality of the support of $\theta_0$, i.e., the set $\{i: \theta_{0i}\ne 0\}$.
	
	\subsection{Recovery using hard-spike-and-slab priors}
	
	It is clear that if each $\theta_i$ is given $\mathrm{Lap}(\lambda)$-prior independently, the resulting posterior mode will be the LASSO, which, with the choice $\lambda\asymp \sqrt{\log n}$, converges to the true $\theta_0$ at the nearly optimal rate $s\log n$. Still, the whole posterior has a suboptimal contraction property (Castillo et al. \cite{csv15}). This is because without sufficient prior concentration at $0$, the posterior cannot contract sufficiently quickly near the truth. A remedy is to assign an additional point-mass at zero using a 
	hard-spike-and-slab prior: for $w\in[0,1]$ and $\Gamma$ a distribution on $\mathbb{R}$,  $\Pi_w = \Pi_{w,\Gamma} = \bigotimes_{i=1}^n\{ (1-w)\delta_0 + w \Gamma\}$. The weight parameter $w$ in the `oracle' situation (i.e., $s$ is known) can be taken to be $s/n$, leading to the optimal rate $s\log (n/s)$. Realistically,  the choice $w=c/n$ is possible with a constant $c>0$, but this leads to a slightly suboptimal rate $s_0\log{n}$, where $s_0$ is the true value of $s$. In this case, under the prior, the expected number of nonzero coefficients is of the order of a constant. To obtain an improved data fit, options that perform better in practice include an empirical Bayes choice $\hat w$ for $w$ (George and Foster \cite{georgefoster}, Johnstone and Silverman \cite{js04}). One possibility is to use an ad-hoc estimator of the number of non-zero coordinates of $\theta$, for instance, by keeping only coordinates above the expected noise level
	$ \hat w =n^{-1} \sum_{i=1}^n \mathbbm{1}\{ |Y_i|> \sqrt{2\log{n}} \}$. 
	However, this choice may be too conservative, as signals below the universal $\sqrt{2\log{n}}$ threshold may go undetected.   
	The marginal maximum likelihood empirical Bayes approach (MMLE) consists of integrating out $\theta$ and maximizing with respect to $w$:  
	$ \hat{w} = \text{argmax}\, \prod_{i=1}^n \left( (1-w)\phi(Y_i) + w g(Y_i) \right)$, where $g=\gamma * \phi$.  
	The plug-in posterior $\Pi_{\hat w}(\cdot|Y_1,\ldots,Y_n)$ was advocated and studied in George and Foster \cite{georgefoster}, Scott and Berger \cite{scottberger10}, Johnstone and Silverman \cite{js04}, Castillo and Mismer \cite{cm18} among others, and is shown to possess the optimal contraction rate for heavy enough tail of the slab.  Castillo and van der Vaart \cite{cv12} considered the hierarchical Bayes using a prior $w\sim \text{Beta}(1,n+1)$. More generally, they also considered 
	a {\em subset-selection prior}: if $\pi_n$ a prior on the set $\{0,\ldots, n\}$ and $\mathcal{S}_s$ the collection of all subsets of $\{1,\ldots,n\}$ of size $s$, let $\Pi$ be constructed as 
	\begin{equation} 
		\label{subsprior}
		s \sim \pi_n, \qquad  S| s  \sim \text{Unif}(\mathcal{S}_s), \qquad 
		\theta| S \sim \bigotimes_{i\in S} \Gamma \,\otimes\, \bigotimes_{i\notin S} \delta_0.
	\end{equation}
	The hard-spike-and-slab prior is a special case where $\pi_n$ is the binomial distribution $\text{Bin}(n,w)$  {(see also van Erven and Szab\'o \cite{szabovanerven20} for further discussion)}.  Through the prior $\pi_n$, it is possible to chose priors on the dimension that `penalize' large dimensions more than the binomial and achieve the optimal rate $s_0\log (n/s_0)$ uniformly over $\ell_0[s_0]$, that is, when the true value of sparsity is $s_0$, for instance, by using the complexity prior $\pi_n(s)\propto \exp[-as\log(bn/s)]$; see Castillo and van der Vaart \cite{cv12}. Finally, deriving the optimal constant $2$ in the minimax sense for posterior contraction is possible for some carefully designed priors as shown by Castillo and Mismer \cite{mismerphd19, cm19}. Song \cite{song20} considered a polynomially-tailed prior with a scaling parameter and showed that, for an appropriate choice of the scaling constant, the posterior contracts at the minimax rate, and the multiplicative constant appearing with the rate is also nearly optimal.  Song \cite{song20} also provided a heuristic explanation of the result by arguing that the Bayesian procedure mimics the hard-thresholding process to choose the signals.

	It may be mentioned that the convergence speeds of the entire posterior distribution and those of aspects of it such as the posterior mean, median or mode may not match. Like the posterior mode for the Laplace prior, the (empirical or hierarchical Bayes) posterior mean for the hard-spike-and-slab prior may also be located `relatively far' from the bulk of the posterior mass in terms of the distance $d_q(\theta,\theta')=\sum_{i=1}^n |\theta_i-\theta'_i|^q$, $0<q<1$ (Johnstone and Silverman \cite{js04}, Castillo and van der Vaart \cite{cv12}). In contrast, the posterior distribution concentrates at the minimax rate for the $d_q$-loss for any $q\in(0,2]$ (Castillo and van der Vaart \cite{cv12}). Further, the coordinate-wise median of the MMLE empirical Bayes posterior for the hard-spike-and-slab prior converges at the minimax rate  (Johnstone and Silverman \cite{js04}). Still, the plug-in posterior $\Pi_{\hat w}(\cdot|Y_1,\ldots,Y_n)$ converges at a suboptimal rate for specific sparse vectors (Castillo and Mismer \cite{cm18}) if a Laplace slab is chosen. This is due to an excessively large dispersion term in the plug-in posterior when $\hat w$ is slightly above $s_0/p$. The problem is avoided in the empirical Bayes approach if a slab with sufficiently heavy tails is chosen (e.g., Cauchy), and in the hierarchical approach through additional regularization via a well-chosen Beta prior on $w$.  
	
	{Computations using hard-spike-and-slab priors in the sequence model are possible in polynomial time. Implementation of the MMLE empirical Bayes approach amounts to a root-finding algorithm and is done in Johnstone and Silverman \cite{js04}; it can handle very large $n$'s (e.g. $n=10^8$). Castillo and van der Vaart \cite{cv12} proved that aspects of the subset-selection posterior induced by \eqref{subsprior}, including hierarchical hard-spike-and-slab priors, are computable in polynomial time. Szab\'o and van Erven \cite{szabovanerven20} provided fast algorithms for $n$ up to the order of $10^5$.}
	
	\subsection{Uncertainty quantification} 
    \label{using normal slab}

If the slab distribution $\Gamma$ is chosen normal, then the posterior is conjugate, conditional on the selection, and hence allows certain explicit expressions, which are beneficial. However, due to quickly decaying tails of normal distributions, this leads to over-shrinkage of large values, leading to suboptimality, or even inconsistency. The problem can be avoided by introducing a mean parameter for the slab distribution of each entry and estimating those by empirical Bayes, which is the same as plugging in the observation itself; see Martin and Walker \cite{martin2014} who used a fractional-posterior, and Belitser and Nurushev \cite{belnur15}, who considered the usual posterior distribution. Belitser and Nurushev \cite{belnur15}  obtained the optimal contraction rate following an oracle approach to optimality, and also obtained credible regions with frequentist coverage and adaptive optimal size under an ``Excessive Bias Restriction''  (EBR) condition that controls bias by a multiple of variability at a parameter value. Restricting parameter space in such a manner is essential as honest coverage with an adaptive-size ball is impossible to obtain by any method, Bayesian or not, in view of impossibility results such as Baraud \cite{baraud02}. Castillo and Szab\'o \cite{CS20} considered the MMLE-empirical Bayes approach using heavy-tailed slabs and an estimated value of the hyperparameter $w$ in the prior for the sparsity level and provided similar coverage results for adaptive credible sets,  obtaining results in terms of $d_q$ distances, $q\le 2$. They also showed that the EBR condition is necessary in a minimax sense. 

\subsection{Alternative shrinkage priors} 

While subset selection priors are particularly appealing because of their natural model selection abilities, to address less strict sparsity (such as weak or strong $\ell_p$-classes), one may use a soft-spike-and-slab prior, to get essentially the same posterior contraction rate, provided that the spike distribution is sufficiently concentrated, such as by the spike-and-slab-LASSO with the spike distribution $\text{Lap}(\lambda_0)$ for $\lambda_0\to \infty$ with $n$ and the slab distribution $\text{Lap}(\lambda_1)$ for constant $\lambda_1$ (Ro\u ckov\`a and George  \cite{rockovageorge18}). For the horseshoe prior, van der Pas et al. \cite{skv14,svv17} obtained explicit expressions for the posterior mean in terms of degenerate hypergeometric functions, and used that to show posterior contraction in terms of the Euclidean distance at any vector in $\ell_0[s_0]$ at the rate $\sqrt{s_0 \log n}$, respectively under known sparsity and unknown sparsity regimes with $0\le s_0\le n$ possibly depending on $n$; see also Ghosh and Chakraborti \cite{ghoshchakra17} for similar results. The optimal posterior contraction rate using the Dirichlet-Laplace prior was obtained by Bhattacharya et al. \cite{bpd15}  under a growth condition 
on the norm $\|\theta_0\|$. The optimal contraction rates for the spike-and-slab-LASSO in its hierarchical form were obtained in Ro\u ckov\`a \cite{rockova18} up to a side-condition on signals, and by Castillo and Mismer \cite{cm18} for the empirical Bayes form. 
Van der Pas et al. \cite{ssv16} unified results for the scale-mixture of normals with a known sparsity parameter, obtaining sufficient and almost necessary conditions for optimal convergence. 
Naturally, to detect the sparsity structure using continuous shrinkage priors, a thresholding procedure that identifies entries essentially zero is needed. Credible $\ell_2$-balls for the parameter and marginal credible with asserted frequentist coverages were obtained in van der Pas et al. \cite{svvuq17}. 

\subsection{Multiple testing}

The problem of identifying a sparse structure can also be viewed as a multiple testing problem, in which case measures such as the false discovery rate (FDR) (Benjamini and Hochberg \cite{BH1995}) may be considered instead of the usually adopted family-wise error rate for Bayesian model selection. The FDR control may fail for many priors. In Castillo and Roquain \cite{CR20}, a uniform FDR control over $\ell_0[s_0]$ is derived for procedures based on both $\ell$-values $\Pi(\theta_i=0|Y_1,\ldots,Y_n)$ and $q$-values $\Pi(\theta_i=0|Y_i\ge y_i)$, when the prior distribution has suitably heavy slab tails and the mixture weight of the spike-and-slab prior is estimated by Empirical Bayes. It is shown that thresholding the $\ell$-values at a given level $t$ leads to the FDR going to $0$ uniformly at a logarithmic rate. The procedure of thresholding $q$-values at a level $t$, on the other hand, controls the FDR at a level a constant times $t$ uniformly, and exactly $t$ asymptotically when signals are all above the threshold $a\sqrt{2\log(n/s_0)}$, for some $a>1$. Abraham et al. \cite{acr20} consider a method based on averages of ordered $\ell$-values as a way to set the rejection threshold. They formally show that it has fairly similar properties as the $q$-value, with an FDR going to the target level $t$ asymptotically for a large class of truths with sufficient signal strength.

The previous multiple testing problem is fully studied from the minimax perspective of control of the sum of FDR and FNR (False Negative Rate, the proportion of false nulls) in the work Abraham et al. \cite{acr24}, where the sharp minimax constant is completely characterised in terms of the sequence of nonzero signals. A surprising take-away from that study is that for sub-Laplace noise (including, in particular, the commonly encountered Gaussian noise setting), for all `reasonable' procedures (in the sense of ones that have several discoveries not significantly larger than $s_0M_n$, for some $M_n\to\infty$), one cannot make the FNR asymptotically strictly smaller even if increasing the FDR significantly. In particular, the Benjamini and Hochberg (BH) procedure (Benjamini and Hochberg \cite{BH1995}) with a {\em fixed} level $t\in(0,1)$ is shown to be suboptimal for the FDR+FNR risk. It is also revealed that the $\ell$--value procedure for a spike-and-slab prior as in the previous paragraph is minimax optimal. In contrast, the BH procedure can achieve this if its level tends to zero, under slightly stronger conditions.  

% {Simultaneous control of FDR and  is considered from a minimax perspective in Castillo and Roquain \cite{CR20} and Abraham et al. \cite{acr20}, where a method based on averages of ordered $\ell$-values as a way to set the rejection threshold is studied. }  
 
Early results for continuous shrinkage priors, including the horseshoe, were obtained in Salomond \cite{salomond17}; see also van der Pas et al. \cite{svvuq17} for a simulation study. The work by Banerjee et al. \cite{bcv26} demonstrates that quantities analogous to the $\ell$--values as defined for spike--and--slab priors above (and that are not available here since the prior puts no atom at $0$) can be obtained from properly chosen posterior marginal quantiles, and are optimal for the FDR+FNR risk in the normal means model; empirical evidence that such procedures also control the FDR is provided also for high-dimensional regression an graphical models.  

 In the Bayesian setting, a very appealing measure is provided by the Bayes risk for testing. With a hard-spike-and-(normal) slab prior with variance $v_1$ and known $w$, the oracle Bayes rule for rejecting $H_{0i}:\theta_i=0$ can be shown to be a thresholding procedure $|Y_i|^2> (1+v_1^{-1})[\log (1+v_1)+ 2\log ((1-w)/w)]$; see Datta and Ghosh \cite{datta2013asymptotic} for details. 

Beyond the sequence model and sparse priors, but in a similar spirit, the Bayes risk for multiple testing in the setting of Hidden Markov Models (HMMs) with two hidden states is studied in Abraham et al. \cite{acg22}. The considered sequence is one of the hidden states, and the $\ell$-values are $\Pi(\theta_i=0|X)$, where $X$ is the observed sequence of emissions from the HMM and $\theta_i=0$ means that the hidden state for observation $i$ has label $0$. The same procedure as in Abraham et al. \cite{acr20} based on cumulative $\ell$-values is shown to be near-asymptotically optimal in terms of true discoveries among procedures controlling the (marginal) Bayesian FDR at level $t$.

\subsection{Multiple change points}
\label{multiple change points}

In some contexts, the distributions of observations vary over time, changing at certain points, and remaining the same over blocks $B_1,\ldots,B_s$ of consecutive locations: $X_1,\ldots,X_n$ are independently distributed as $F_1,\ldots,F_n$ 
respectively, with an unknown number $s-1$ of change points $1< i_1<\cdots <i_{s-1}< n$. 
When $s\ll n$, estimating $s$ and $i_1,\ldots,i_s$ is possible. A common modeling assumption is that $F_i=\mathrm{N}(\theta_i,\sigma^2)$, $i=1,\ldots,n$, for some common but possibly unknown $\sigma$. Then, studying the multiple change-point model reduces to making an inference on a long normal mean vector. However, rather than sparsity, the hidden lower-dimensional structure is the constancy of the mean across blocks, which is equivalent to the sparsity of the successive differences vector. In the classical literature, approaches based on multiple testing, various 
CUSUM versions, (variants of) binary segmentation, multiple testing with varying windows, step-wise regression, 
screening and ranking algorithms, multi-scale methods. More recently, high-dimensional sparse regression representation, exploiting the sparsity of $\theta_i-\theta_{i-1}$ through the fused lasso penalty $\lambda \sum_{i=2}^n |\theta_i-\theta_{i-1}|$ has been popular; see Tibshirani et al. \cite{tibshirani2005sparsity}, Rinaldo \cite{rinaldo2009properties},
and the review articles Carlstein et al. \cite{carlstein1994change}, Chen and Gupta \cite{chen2012parametric}, Lee \cite{lee2010change}, and Niu et al. \cite{niu2016multiple}. 

An empirical Bayes method was proposed by Belitser and Ghosal \cite{belitser2025bayesian}. 
Apart from the posterior contraction rate and the accuracy of change-point selection under a minimum separation condition for the means in the successive blocks, they also provided results on adequate coverage, with size adapting to $s$. They placed a prior on $s$ exponential in the model complexity measure $s\log(n/s)$ and a conjugate normal prior with block-means estimated from the data, and obtained an oracle bound for posterior contraction in terms of the oracle risk at a given parameter vector to derive these properties, using an oracle inequality of Belitser and Nurushev \cite{belitsernur19} based on the approach in Belitser \cite{belitser17}. The approach does not need to assume that the actual data-generating mechanism is normal with a known variance; it only requires a sub-Gaussianity condition with a known bound on the variance. In particular, the observations can be dependent. The size of the credible ball with adequate coverage adapts to the optimal order of the size $\sqrt{s\log (n/s)}$ under the EBR condition, as in the case of a sparse sequence; here, as in Subsection~\ref{using normal slab}, the EBR condition controls the bias by a multiple of the variability at a parameter value, but relative to the block structure rather than the sparse structure. To detect the change points, they showed that if $s\lesssim n^{\delta}$ for some $\delta<1$, then the change points are correctly identified with probability approaching 1 at a polynomial rate in $n^{-1}$ if the minimum separation of the means in successive blocks are at least a sufficiently large multiple of $\sigma \sqrt{1+\log n}$, and the rate of separation cannot be improved without additional assumptions. However, if the block sizes are assumed to be at least $L$ and the number of blocks is $O(\log n)$, then the rate of separation can be improved by a factor of $L^{-1/2}$, provided the zero-one loss of change points identification is replaced by a Hausdorff-type loss function. 

A fully Bayesian alternative to the $\ell_1$-based fusion approach is to induce sparsity directly on the successive differences $\eta_i = \theta_i - \theta_{i-1}$ via shrinkage priors. Early formulations used Laplace or heavy-tailed priors, such as the $t$-shrinkage prior of Song and Cheng \cite{song2020bayesian}, which improve robustness over $\ell_1$-penalties while still imposing relatively uniform shrinkage. A more flexible approach is the horseshoe fusion prior of Banerjee \cite{banerjee2022horseshoe}, which places a global-local horseshoe prior on $\eta_i$, enabling strong shrinkage of small differences while preserving large jumps corresponding to true change points. Owing to its unbounded spike at zero and heavy tails (Carvalho et al. \cite{carvalho}, Polson and Scott \cite{Polson&Scott}), the horseshoe prior reduces bias and adapts effectively to the unknown block structure. The resulting Bayesian formulation provides full posterior inference for the mean sequence and change points, and achieves posterior contraction and consistent recovery under suitable separation conditions. Extensions of this framework to graph-structured signals, where fusion is defined over general graphs rather than linear chains, have also been developed; see, for instance, Banerjee and Shen \cite{banerjee2022graph} and Banerjee \cite{banerjee2022horseshoe}.

\section{High-dimensional regression}
\label{regression}

A natural generalization of the normal sequence model is the situation where 
the mean $\theta_i$ depends on a (high-dimensional) covariate $X_i$ associated with the observation $Y_i$, $i,\ldots,n$.
The most popular statistical model in this setting is given by the normal linear regression model $Y_i=\beta^{\mathrm{T}} X_i+\varepsilon_i$, $i=1,\ldots,p$, where $X_i\in \mathbb{R}^p$ is a deterministic $p$-dimensional predictor, $\beta \in \mathbb{R}^p$ is the linear regression coefficient, and $\varepsilon_i$ are i.i.d. $\mathrm{N}(0,\sigma^2)$ error variables. We often write the linear regression model in the equivalent vector form 
\begin{align}
    \label{linear model} 
    Y=X\beta+\varepsilon, \quad Y,\varepsilon\in \mathbb{R}^n, \; X\in \mathbb{R}^{n\times p},\; \varepsilon\sim \mathrm{N}_n(0, \sigma^2 I_{n\times n}).
\end{align}

When $(p^4 \log p)/n\to 0$, posterior concentration and a Bernstein-von Mises theorem were obtained in Ghosal \cite{ghosal1999asymptotic} without assuming any sparsity. In the high-dimensional setting where $p$ is very large (possibly even much larger than $n$), the behavior of the posterior distribution of $\beta$ has been studied assuming that $\beta$ is sparse with only $s$ coordinates of it being non-zero, $s\ll n$, and $n\to \infty$. 

\subsection{Linear regression with hard-spike-and-Laplace slab} 

Let $X_{.,j}$ be the $j$th column of the matrix $X:=(\!( X_{ij})\!)$, and 
consider the norm 
$\|X\|  =\max \|X_{.,j}\|$.  
We consider the prior \eqref{subsprior}, with $g_S$ the product of $|S|$ Laplace densities $\beta\mapsto (\lambda/2) \exp(-\lambda|\beta|)$. We set $\lambda=\mu\|X\|$, for 
$p^{-1} \le \mu \le 2\sqrt{\log p}$ and put 
a prior $ \pi_p(s) \propto c^{-s}p^{-as}$ on the number of non-zero entries of $\beta$. For $p>n$, clearly $\beta$ cannot be uniquely recovered from $X \beta$ (even in the noiseless case). However, if $\beta$ is assumed to be sparse (with only $s\ll n$ components non-zero) and the submatrix of $X$ corresponding to the active predictors is full rank $s$, it turns out that $\beta$ can be recovered. Define the \emph{compatibility number of model} $S\subset \{1,\ldots,p\}$ by 
$$
\overline\phi(S):=\inf\Big\{{\|X\beta\| |S|^{1/2}}/{\|X\|\, \|\beta_S\|_1 }
: \|\beta_{S^c}\|_1\leq 7 \|\beta_S\|_1,  \, \beta_S\neq 0\Big\}, $$
where $\beta_S=(\beta_i: {i\in S})$, 
and the \emph{$\ell_r$-compatibility number} in vectors of dimension $s$ by 
$$ \phi_r(s):=\inf\Big\{{\|X\beta\| |S_\beta|^{1-r/2}}/{ \|X\|\, \|\beta\|_r }:\,
0\neq |S_\beta| \leq s\Big\}, \quad r=1,2,$$ 
where $S_\beta=\{j: \beta_j\ne 0 \}$, the support of a sparse vector $\beta$. For $\beta_0$ the true vector of regression coefficients and $S_0=S_{\beta_0}$ and $s_0=|S_0|$, we assume that, with $r=1$ or $2$ depending on the context,  
$\min\left(\phi_r(S_0),\overline\phi(Cs_0)\right)\ge d>0$, 
where $C$ is a suitably large constant depending on $S_0,a,\mu$. More information about compatibility numbers may be found in van de Geer and B\"uhlmann \cite{van2009conditions}. When the entries of $X$ are sampled randomly, the compatibility numbers are often very well-behaved with high probability. In the leading example where $X_{ij}$ are independently sampled from  $\mathrm{N}(0,1)$,  the compatibility numbers above of models up to the dimension a multiple of $\sqrt {n/\log p}$ are bounded away from zero (Cai and Jiang \cite{caijiang11}, van de Geer and Muro \cite{geermuro14}).

The following conclusions were derived in Castillo et al. \cite{csv15}. 
First, a useful property is that the posterior of $\beta$ sits on models of dimensionality at most a constant multiple of the true dimension $s_0$, thereby returning models at least of the same order of sparsity as the true one. Further, we can recover $\beta_0$ in terms of the $\ell_1$-norm: when the $\ell_1$-compatibility numbers are bounded away from $0$,   
$$\mathrm{E}_{\beta_0}\Pi ( \|\beta-\beta_0\|_1 > M s_0  \sqrt{\log  p}/{\|X\|}|Y_1,\ldots,Y_n)  \to 0. $$
The corresponding rate in terms of the Euclidean distance is $\sqrt{s_0\log  p}/{\|X\|}$, assuming that the $\ell_2$-compatibility numbers are bounded away from $0$. The corresponding rate is $\sqrt{\log  p}/{\|X\|}$ with respect to the maximum-loss under a stronger condition, called \emph{mutual coherence}.
The rate for prediction, i.e., bound for $\|X\beta-X\beta_0\|$, is of the order $\sqrt{s_{0}  \log  p}$ under a slightly adapted compatibility condition. The convergence results are uniform over the parameter space under boundedness conditions on the compatibility numbers and match those of celebrated estimators in the frequentist literature.    

In sparse regression, variable selection—the identification of non-zero coefficients—is essential because it enables simpler interpretation and a better understanding of relations. To address the issue, Castillo et al. \cite{csv15} developed a technique based on a distributional approximation under relatively low choices of the parameter $\lambda$, known as the small-$\lambda$ regime. This is similar to the normal approximation to the posterior distribution in the Bernstein-von Mises theorem. Still, the difference is that in this context, the approximating distribution is a mixture of sparse normal distributions over different dimensions. Then the problem of variable selection can be transferred to that for the approximate posterior distribution, and it can be shown that no proper superset of the correct model can be selected with appreciable posterior probability. It is, however, possible to miss signals that are too small in magnitude. If all non-zero signals are assumed to be at least of the order $\sqrt{s_0\log p}/\|X\|$, then none of these signals can be missed, because missing any of them introduces an error of the magnitude of the contraction rate. Thus, the selection consistency follows. Also, in this situation, the distributional approximation reduces to a single normal component, with sparsity exactly as in the true coefficient. Ning et al. \cite{ning2018bayesian} considered a useful extension of the setting of Castillo et al. \cite{csv15} by letting each response variable $Y$ to be $d$-dimensional with $d$ also increasing to infinity (at a sub-polynomial growth with respect to $n$), and having completely unknown $d\times d$-covariance matrix $\Sigma$ and the regression has group-sparsity. They used a prior on $\Sigma$ using a Cholesky decomposition and used a hard-spike-and-slab prior with a multivariate Laplace slab on the group of regression coefficients selected together. Applying the general theory of posterior contraction (Ghosal et al. \cite{ggv}) using exponentially consistent tests for separation based on the R\'enyi divergence, they obtained the squared posterior contraction rate 
\begin{equation}
	\label{rate group sparsity}
	\epsilon_n^2=\max \left\{ (s_0 \log G)/n, (s_0 p_{\max}\log n)/n, (d^2 \log n)/n\right\}\!,
\end{equation}
where $G$ stands for the total number of groups of predictors, $s_0$ the number of active groups, and $p_{\max}$ the maximum number of predictors in a group, provided that the regression coefficients and the covariance matrix are appropriately norm-bounded. 
Ning et al. \cite{ning2018bayesian} further extended the distributional approximation technique of Castillo et al. \cite{csv15} to demonstrate selection consistency. 
A general setup of sparse linear regression 
$$Y_i=X_i \beta+\xi_{\eta,i}+\varepsilon_i, \quad \varepsilon_i\sim \mathrm{N}_{m_i}(0,\Delta_{\eta,i})$$ 
independently with possibly varying dimension $m_i$ and covariance matrices $\Delta_{\eta,i}$, allowing a  nuisance parameter $\eta$ and an additional term $\xi_{\eta,i}$ to incorporate various departures from the simple linear model $Y_i=X_i\beta+\varepsilon_i$, was considered in Jeong and Ghosal \cite{JeongLinear}. They showed optimal recovery for the regression coefficient $\beta$ under standard conditions on compatibility numbers. They established recovery rates from the posterior contraction rate for both the regression coefficient $\beta$ in terms of the Euclidean norm and for the nuisance parameters in terms of the squared distance $n^{-1} \sum_{i=1}^n \{\|\xi_{\eta,i}-\xi_{\eta',i}\|^2+ \|\Delta_{\eta,i}-\Delta_{\eta',i}\|_{\mathrm{F}}^2\}$, which they derived by applying the general theory of posterior contraction in terms of the R\'enyi divergence as in Ning et al. \cite{ning2018bayesian}. Further, extending the distributional approximation technique, Jeong and Ghosal \cite{JeongLinear} also showed selection consistency. Their setup accommodates a variety of practical extensions of the basic linear model, including multidimensional response, partially sparse regression, multiple responses with missing observations, (multivariate) measurement errors, parametric correlation structure, mixed-effect models, graphical structure in precision matrix (see Subsection~\ref{subsec:precision matrix}), nonparametric heteroscedastic regression, and partial linear models. 

{While the previous framework enables one to obtain asymptotic normality for {\em all} coordinates in $S_0$ if these are suitably large,  if one is interested in just one coordinate say $\beta_1$, it is possible to build an estimator (or marginal posterior) thereof that is asymptotically normal and efficient, under say only a compatibility--type condition: this has been obtained by Zhang and Zhang \cite{zz14} for a debiased--LASSO estimator in the regime $s\log{p}=o(\sqrt{n})$. In Yang \cite{yang19}, the author derives an analogous Bayesian asymptotic normality result for the marginal posterior $\beta_1$ given $Y$. This is achieved by first defining a suitable prior distribution on $\beta_1$ given the remaining coordinates $\beta_{-1}$, and then assigning a parsimonious prior to $\beta_{-1}$ as before.
	
	%Suppose one is interested in $\theta_1$ and denote by $\theta_{-1}=(\theta_2,\ldots,\theta_p)$. Let $\gamma_i=X_1^t X_1/\|X_1\|^2$ and $W=(I-H)X$, where $H$ is the projection matrix onto $\text{Span}(X_1)$. 
	%The idea is to put a prior with a certain dependence structure on $(\theta_1,\theta_{-1})$ by setting 
	%\[ \theta_1| \theta_{-1}\sim N\Big(-\sum_{i=2}^p \gamma_i\theta_i,\sigma_n^2\Big),\qquad
	%\theta_{-1} \sim \Pi_W,
	%\]
	%where $\sigma_n^2$ %\gg\|\te\|_1\la_n/\|X_1\|_2$
	%is suitably large and $\Pi_W$ a sparsity prior (such as the one of \cite{csv15} or \cite{gvz15}) built as if the input design matrix would be $W$ instead of $X$. Then, under similar conditions as in \cite{zz14}, for $d_{BL}$ the bounded Lipschitz metric,
	%\[ d_{BL}(\|X_1\|_2(\theta_1-\hat\theta_1),N(0,1))\to^P 0,\]
	%where $\hat\theta_1$ is an efficient estimator of $\theta_1$.
}

\subsection{Linear regression using other priors}

Instead of the hard-spike-and-slab prior with a Laplace slab, we may use a normal slab with an undetermined mean selected via empirical Bayes, as in the sequence model. The advantage is that posterior conjugacy, given the chosen subset, can be retained, allowing explicit expressions. Belitser and Ghosal \cite{belitserghosal16} followed the oracle approach of Belitser \cite{belitser17} to quantify risk locally at each point, thereby defining the optimal concentration, and established results analogous to those of Castillo et al. \cite{csv15}. Moreover, uniformly over the collection of all parameter values satisfying the excessive-bias restriction, they showed that an appropriately inflated Bayesian credible ball of the optimal size has adequate frequentist coverage for any sparsity level. Other references on Bayesian high-dimensional linear regression include Narisetty and He  \cite{narisettyhe14}, who considered an alternative to the posterior measure called the {\em skinny Gibbs posterior}, which avoids a critical computational bottleneck and establishes contraction at the true coefficient at the optimal rate.  

An alternative to using an independent hard-spike-and-Laplace slab is to use an elliptical Laplace prior on the selected coefficients through a set selection prior in the spirit of \eqref{subsprior}, but adjusted for the normalizing constant appearing in the elliptical Laplace distribution. Gao et al. \cite{gvz15} considered this prior for a structured linear model $Y=L_X \beta+\varepsilon$, where $L_X$ is a linear operator depending on $X$ and $\varepsilon$ is a vector of errors assumed to have only sub-Gaussian tails, and obtained minimax contraction rates. Apart from linear regression, their setup includes stochastic block models, biclustering, group-sparse linear regression, multi-task learning, dictionary learning, and others. 

Despite the attractive theoretical properties, posterior computation based on two-component priors is computationally intensive. Faster computation may be possible with continuous shrinkage priors (Bhattacharya et al. \cite{bhattacharya2016fast}). Under certain general conditions on prior concentration near zero, the thickness of the tails and additional conditions on the eigenvalues of the design matrix, Song and Liang \cite{song2017nearly} derived posterior contraction and variable selection properties. Their results cover a wide range of continuous shrinkage priors, including the horseshoe, Dirichlet-Laplace, normal-gamma, and t-mixtures, and provide optimal contraction rates and selection consistency (modulo an appropriate thresholding operation) under appropriate conditions on the shrinkage parameters. 

\subsection{Regression beyond linear} 

Regression beyond the linear setting in the high-dimensional context is of substantial interest. The most common extension is a generalized linear model (GLM). A conjugate prior for Bayesian inference in the GLM setting was introduced by Chen and Ibrahim  \cite{chen2003conjugate}. The prior was further used for variable selection in a high-dimensional GLM by Chen et al. \cite{chen2008bayesian}. One of the first papers on convergence properties of posterior distributions in the high-dimensional Bayesian setting is Jiang \cite{jiang2007bayesian}, who derived consistency of the posterior distribution in a GLM with the dimension $p$ possibly much larger than the sample size $n$, needing only $\log p=o(n^{1-\xi})$ for some $\xi>0$, but only in terms of the Hellinger distance on the underlying densities. Posterior concentration and a Bernstein--von Mises theorem for the parameter were obtained earlier by Ghosal \cite{ghosal1997normal} without assuming any sparsity structure but under the restricted growth condition $(p^4 \log p)/n\to 0$.  Atchad\'e \cite{atchade2017contraction} considered a pseudo-posterior distribution in a general setting, assuming only a local expansion of the pseudo-likelihood ratio and derived posterior contraction rates for hard-spike-and-slab priors by constructing certain test functions. Jeong and Ghosal \cite{JeongGLM} obtained recovery rates for regression coefficients in a sparse generalized linear model using hard-spike-and-slab priors. For the special case of logistic regression, Wei and Ghosal \cite{wei2017contraction} established posterior contraction utilizing a variety of continuous shrinkage priors using Atchad\'e's test construction and prior concentration bounds of Song and Liang \cite{song2017nearly}. Both posterior contraction rates and selection consistency were established in a more general multidimensional logistic regression model in Jeong \cite{JeongLogistic} using hard spike-and-slab priors. 

Bayesian high-dimensional regression in the nonparametric setting has been addressed for an additive structure with random covariates. Yang and Tokdar \cite{yangtokdar15} considered Gaussian process priors for each selected component function aided by a variable selection prior on the set of active predictors and showed that the minimax rate $\max \{ \sqrt{(s_0 \log p)/n}, \sqrt{s_0} n^{-\alpha/(2\alpha+1)}\}$ is obtained up to a logarithmic factor, where $s_0$ is the number of active predictors and $\alpha$ is the smoothness level of each component function; see also Suzuki \cite{suzuki12} for earlier results in a PAC-Bayesian setting, where `PAC' stands for `Probably Approximately Correct', and Castillo and Randrianarisoa \cite{cr25} who showed that variable selection can also be performed automatically through priors on the Gaussian process lengthscale parameters, without the need for an explicit model selection prior.  Using an orthogonal basis expansion technique, Belitser and Ghosal \cite{belitserghosal16} extended their oracle technique from linear to additive nonparametric regression and obtained the minimax rate with a hard-spike-and-normal slab prior, with the means of the 
slabs selected by the empirical Bayes technique. They also obtained coverage for adaptive-size balls of functions under an analog of the {\it excessive bias restriction} (EBR) condition for random covariates, called the $\epsilon$-EBR condition. Wei et al. \cite{wei2018sparse}, extending the work of Song and Liang \cite{song2017nearly}, obtained the optimal contraction rate and consistency of variable selection for additive nonparametric regression using a B-spline basis expansion prior with a multivariate version of the Dirichlet-Laplace continuous shrinkage prior on the coefficients. The multivariate version is needed to maintain the group-selection structure of the coefficient vector corresponding to the same component function. 

A high-dimensional single-index nonparametric regression model where the regression function is represented as $a_0(X^{\mathrm{T}} \beta, Z^{\mathrm{T}} \eta)-F(t) a_0(X^{\mathrm{T}} \beta, Z^{\mathrm{T}} \eta)$, with $X$ a high-dimensional predictor, $Z$ a low-dimensional predictor and $a_0$, $a_1$, $F$ smooth functions, was considered by Roy et al. \cite{Arka} to model the atrophy of different brain regions over time. For identifiability of the model, the coefficients $\beta$ and $\eta$ need to be unit vectors in respective dimensions, and $\beta$ should be sparse as well in an appropriate sense. To address the problem, Roy et al. \cite{Arka} introduced a notion of sparsity in the polar coordinate system using a soft spike-and-slab prior with a uniform slab and spike distributions with spikes at appropriate multiples of $\pi/2$. They characterized the posterior contraction rate in terms of the average squared distance on the functions using the general theory of posterior contraction essentially as the weakest of the rates $n^{-\iota/(2\iota+2)}$, $n^{-\iota'/(2\iota'+2)}$ and $\sqrt{(s_0\log p)/n}$, up to a logarithmic factor, where $\iota$ is the smoothness of the functions $(a_0,a_1)$, $\iota'$ is the smoothness of the function $F$, and $s$ is the sparsity of the true $\beta$.  

Shen and Ghosal \cite{shen&ghosal} considered the problem of estimating the conditional density of a response variable $Y\in [0,1]$ on a predictor $X\in [0,1]^p$ with $\log p=O(n^{\alpha})$ for some $\alpha<1$, assuming that only $s_0$ predictors are active, where $s_0$ is unknown but does not grow. They put a prior on the conditional density $p(y|x)$ through tensor products of B-spline expansions. First, a subset selection prior with sufficient control over the effective dimension is used to select the active set, and then priors on the lengths of the spline bases are placed. Finally, the coefficient arrays are given independent Dirichlet distributions of the appropriate dimensions. Shen and Ghosal \cite{shen&ghosal} showed the double adaptation property --- the posterior contraction rate adapts to both the underlying dimension of the active predictor set and the smoothness of each function $\alpha$, i.e., $n^{-\alpha/(2\alpha+s_0+1)}$. In fact, the conditional density can have anisotropic smoothness (i.e., different smoothness in different directions), and then the harmonic mean $\alpha^*$ will replace $\alpha$. Norets and Pati \cite{norets2017adaptive} obtained analogous results for the conditional density of $Y\in \mathbb{R}$ on a predictor $X\in \mathbb{R}^p$ using Dirichlet process mixtures of normal densities.

\section{Learning structural relationship among many variables}
\label{sec:relations}

In this section, we review methods for identifying relationships among a large number of variables by examining features of their joint distributions and the associated convergence results. In Subsection~5.1, large covariance matrices are considered. In Subsection~5.2, Bayesian estimation of sparse precision matrices, which may have specific structures like banding or can be unstructured, is considered. Classification using a high-dimensional predictor is considered in Subsection~5.3. Relations among non-Gaussian variables are studied in Subsection~5.4, and semiparametric models for finding structural relations are considered in Subsection~5.5. 

\subsection{Estimating large covariance matrices}

Understanding the dependence among a large collection of variables $X=(X_i: i=1,\ldots,p)$ is an important problem to study. The simplest measure of dependence is the pairwise covariance between these variables, which leads to the problem of estimating a large covariance matrix. In the setting of a large collection of variables, it is natural to postulate that most of these variables are pairwise independent, or at least uncorrelated, meaning that most off-diagonal entries are zero, or at least close to zero. This introduces a useful sparsity structure that allows meaningful inference with relatively fewer observations. One particularly prominent structure is (approximate) banding, which means that when the variables are arranged in some natural order, the pairwise covariances are zero (or decay quickly) if the lag between the corresponding indices exceeds some threshold. For example, an exact banding structure arises in a moving-average (MA) process. In an autoregressive (AR) or autoregressive moving average (ARMA) process, the covariances decay exponentially with lag. Banding or tapering of the sample covariance is often used to estimate such a large covariance matrix. When the sparse structure lacks a specific pattern, threshold methods are often used, but positive definiteness may not be preserved. Another important low-dimensional structure is given by a sparse plus low-rank decomposition of a matrix $\Sigma =D + \Lambda \Lambda^{\mathrm{T}}$, where $D$ is a scalar matrix and $\Lambda$ is a ``thin matrix'', that is, $\Lambda$ is $p\times r$, where $r\ll p$. Moreover, the columns of $\Lambda$ themselves may be sparse, allowing further dimension reduction. Such a structure arises in structural linear models $X_i =\Lambda \eta_i +\varepsilon_i$, where $\Lambda$ is a $p\times r$ sparse factor loading matrix and $\eta_i\sim \mathrm{N}_p(0,I)$ are the latent factors independent of the error sequence $\varepsilon_i \sim \mathrm{N}_p(0,\sigma^2 I)$. In this setting, a Bayesian method was proposed by Pati et al. \cite{pati2014posterior}. They considered an independent hard-spike-and-slab prior with a normal slab on entries of $\Lambda$ and an inverse-gamma prior on $\sigma^2$, along with a Poisson-tailed prior on the rank $r$. They showed that if the number of entries in each column of $\Lambda$ is bounded by $s_0$, and the true $\sigma^2$ and the number of latent factors are bounded, then the posterior contraction rate $\sqrt{(s_0 \log n \log p)/n}$ can be obtained. In a recent work involving covariance matrices with an arbitrary sparsity pattern, Lee et al. \cite{lee2022beta} proposed a beta-mixture shrinkage prior that is more efficient than the spike-and-slab prior and also attains minimax optimality in high dimensions.

\subsection{Estimating large precision matrices and graphical models}
\label{subsec:precision matrix}

An intrinsic relation among a set of variables is described by the conditional dependence of a pair when the effects of the remaining variables are eliminated by conditioning on them. In a large collection of variables, most pairs are conditionally independent given the others. It is convenient to describe this structure using a graph, where each variable corresponds to a node, and an edge connects two nodes if and only if the two are conditionally dependent. Therefore, such models are popularly called graphical models. An introduction to graphical models is given in the Appendix. If $X_i$ and $X_j$ are conditionally independent given $X_{-i,-j}:=(X_k: k\ne i,j)$, then it follows that $\omega_{ij}=0$, where $(\!( \omega_{ij})\!)=\Sigma^{-1}$ is the precision matrix of $X$, to be denoted by $\Omega$. In a Gaussian graphical model (GGM), i.e., when $X$ is distributed as jointly normal, $X\sim \mathrm{N}_p (0, \Omega^{-1})$, then the converse also holds, namely $\omega_{ij}=0$ implies that $X_i$ and $X_j$ are conditionally independent given $X_{-i,-j}$. Thus, in a Gaussian graphical model, the problem of learning the intrinsic dependence structure reduces to estimating the precision matrix $\Omega$ under sparsity (i.e., most off-diagonal entries of $\Omega$ are $0$). 

\subsubsection{Banding and other special sparsity patterns}

A special type of sparsity is given by an approximate banding structure of $\Omega$. Note that this is different from banding of the covariance matrix $\Sigma$, as the inverse of a banded matrix is only approximately banded. Among familiar time series, an AR process has a banded precision matrix, whereas an MA process has an approximate banded structure. The graph corresponding to a banded precision matrix has edge set $\{(i,j): |i-j|\le k, i\ne j\}$, $k$ is the size of the band, which always corresponds to a decomposable graph, allowing the use of structurally rich priors in this context. In particular, a conjugate graphical Wishart (G-Wishart) prior (see the Appendix for its definition and properties) can be put on $\Omega$ together with a choice of $k$, which may even be given a prior. Banerjee and Ghosal \cite{banerjee2014posterior} showed that with this prior, the posterior contraction rate in terms of the spectral norm at an approximately banded true $\Omega$ with eigenvalues bounded away from $0$ and infinity is given by $\max \{k^{5/2} \sqrt{(\log p)/n}, k^{3/2}\gamma(k)\}$, where $\gamma(k)$ is the banding approximation rate given by the total contribution of all elements outside the band. This, in particular, implies that for a $k_0$-banded true precision matrix with fixed $k_0$, the rate is nearly parametric $\sqrt{(\log p)/n}$. The rate calculation extends to general decomposable graphs, with $k$ standing for the maximal cardinality of a clique, as shown by Xiang et al. \cite{xiang2015high}, who used distributional results on the Cholesky factor of decomposable G-Wishart matrices. Another proof of the same result using the original technique of Banerjee and Ghosal \cite{banerjee2014posterior} was given by Banerjee  \cite{banerjee2017posterior}. 

Lee and Lee \cite{lee2017estimating} considered a class of bandable precision matrices in the high-dimensional setup, using a modified Cholesky decomposition (MCD) approach $\Omega = (I_p - A)^{\mathrm{T}} D^{-1}(I_p - A)$, where $D$ is a diagonal matrix, and $A$ is lower-triangular with zero diagonal entries. This results in the {\em $k$-banded Cholesky prior} with non-zero entries of $A$ getting the improper uniform prior and the diagonal entries a truncated polynomial prior. They adopted a decision-theoretic framework and demonstrated convergence in posterior expected loss. Their rate is sharper with $k^{5/2}$ replaced by $k^{3/2}$, but the results are not exactly comparable, as their classes of true precision matrices are different. 
The two become comparable only for near-bandable matrices with exponentially decreasing decay functions, in which case the rates are essentially equivalent. 
Lee and Lin \cite{lee2018bayesian} proposed a prior distribution that is tailored to estimate the bandwidth of large bandable precision matrices. They established strong model-selection consistency for the bandwidth parameter, along with consistency of Bayes factors.

A natural fully Bayesian approach to graphical structure selection is to put a prior $p(G)$ on the underlying graph $G$ and then a G-Wishart prior $p(\Omega|G)$ on the precision matrix $\Omega$ given $G$. 
The joint posterior distribution of $(\Omega,G)$ is then given by
$$
p(\Omega,G|X^{(n)}) \propto p({X^{(n)}}|\Omega,G)p(\Omega| G)p(G),$$ 
where $X^{(n)}$ stand for $n$ i.i.d. observations on $X$. 
Taking a discrete uniform distribution over the space $\mathcal{G}$ of decomposable graphs for $p$ variables, we get the joint distribution of the data and the graph $G$, after integrating out $\Omega$ as
$$
p({X}^{(n)}, G) = ({(2\pi)^{np/2}\#\mathcal{G}})^{-1}{I_G(\delta + n,D+nS_n)}/{I_G(\delta,D)},
$$
where $I_G$ is defined in \eqref{G-Wishart normalization}. This immediately yields an expression for the Bayes factor on which model selection may be based. Computational techniques in graphical models are discussed in the Appendix. 

\subsubsection{Models without specific sparsity patterns}
\label{wo specific sparisty precision}

When the graph lacks an orderly pattern, the precision matrix may be estimated by imposing sparsity on its off-diagonal entries. The most commonly known non-Bayesian method is the \emph{graphical LASSO}  (Friedman et al. \cite{friedman2008sparse}), which is obtained by maximizing the log-likelihood subject to the $\ell_1$-penalty on the off-diagonal elements. A Bayesian analog, called the \emph{Bayesian graphical LASSO}, was proposed by Wang \cite{wang2012bayesian} by imposing independent Laplace priors on the off-diagonal elements and exponential priors on the diagonal elements, subject to a positive definiteness restriction on $\Omega$. As in the Bayesian LASSO, a block Gibbs sampling method based on the Laplace distribution's representation as a scale mixture of normals is used. Clearly, the Bayesian graphical LASSO is motivated by the desire to make the posterior mode the graphical LASSO. The main drawback is that the whole posterior is never supported on sparse precision matrices, and the posterior concentration is suboptimal, as in the case of the Bayesian LASSO. The forceful restriction to positive definiteness also leads to an intractable normalizing constant in the prior. However, Wang \cite{wang2012bayesian,wang2015scaling} developed a clever computational trick to avoid the computation of the normalizing constant in posterior sampling, known as \emph{scaling-it-up}. More details of the computational algorithm are described in the Appendix. To alleviate the problem of the lack of sparsity in the posterior, Wang \cite{wang2012bayesian} used a thresholding approach proposed in Carvalho et al. \cite{carvalho}. Li et al. \cite{li2019graphical} proposed using the horseshoe prior on the off-diagonal elements instead of Laplace, thus leading to the `graphical horseshoe' procedure. They also developed an analogous block Gibbs sampling scheme using a variable-augmentation technique for half-Cauchy priors, as proposed in Makalic and Schmidt \cite{makalic2015simple}. The resulting procedure appears to have a better posterior concentration, as measured by the Kullback-Leibler divergence, and smaller bias for non-zero elements. Wang and Pillai \cite{wang2013class} explored other shrinkage priors based on a scale-mixture of uniform distributions. In a recent work, Sagar et al. \cite{sagar2024precision} developed a graphical horseshoe-like prior for fast estimation of precision matrices using a novel prior-penalty dual approach.

The convergence properties of the posterior distribution of a sparse precision matrix with an arbitrary graphical structure were studied by  Banerjee and Ghosal \cite{banerjee2015bayesian}. They considered a prior similar to the Bayesian graphical LASSO except that a large point-mass at $0$ is added to the prior distribution of off-diagonal entries, and the total number of non-zero off-diagonal entries is restricted. They derived posterior contraction rate in terms of the Frobenius norm as $\sqrt{((p+s_0)\log p)/n}$, where $s_0$ is the number of true non-zero off-diagonal entries of the precision matrix, which is the same as the convergence rate of the graphical LASSO and is optimal in that class. The proof uses the general theory of posterior contraction (Ghosal et al. \cite{ggv}) by estimating the prior concentration near the truth and the entropy of a suitable subset of precision matrices, which receives most of the prior mass. The sparsity built into the prior helps control the effective dimension, and thus the prior concentration and entropy.  Control over the eigenvalues allows linking the Frobenius norm to the Hellinger distance for normal densities. Their technique of proof extends to soft-spike-and-slab and continuous shrinkage priors,  as long as the prior concentration at zero off-diagonal values is sufficiently high; for example, Sagar et al. \cite{sagar2024precision} established optimal posterior contraction rates for both the graphical horseshoe and the graphical horseshoe-like priors. In a Gaussian graphical model with observations subject to independent normal measurement error, Shi et al. \cite{shi2021bayesian} showed that the rate obtained in Banerjee and Ghosal \cite{banerjee2015bayesian} prevails even if the measurement error is incorporated in the model, provided the true measurement error variance is below a particular threshold. They considered several prior distributions for both spike-and-slab and continuous shrinkage prior types, either directly on the entries or via a sparse Cholesky factorization, to different kinds of true sparsity structures, such as general sparsity, banding, or a factor model. They also introduced a key technique that adds an independent inverse-Gaussian nugget factor to the diagonal elements. This avoids a condition that the prior distribution of the minimum eigenvalue of the precision matrix is floored, thus allowing the use of many common sparse priors for the precision matrix without an artificial truncation step. They also showed that ignoring the measurement error leads to inconsistency. 

Niu et al. \cite{niu2019bayesian} addressed the problem of graph selection consistency under model misspecification. In the well-specified case where the true graph is decomposable, strong graph selection consistency holds under certain assumptions on the graph size and edge counts, using a G-Wishart prior. For the misspecified case, where the true graph is non-decomposable, they showed that the posterior distribution of the graph concentrates on the set of minimum triangulations of the true graph.

\subsection{High-dimensional discriminant analysis}
\label{discriminant}

Consider a classification problem based on a high-dimensional  predictor $X=(X_1,\ldots,X_p)\sim \mathrm{N}_p(\mu,\Omega^{-1})$, where $(\mu,\Omega)=(\mu_1,\Omega_1)$ for an observation from the first group and $(\mu,\Omega)=(\mu_2,\Omega_2)$ when the observation comes from the second group. Linear and quadratic discriminant analyses are popular model-based classification methods. The oracle Bayes classifier, which uses the true values of $\mu_1,\mu_2$ and $\Omega_1,\Omega_2$, has the lowest possible misclassification error. Using a Bayesian procedure with priors on $\mu_1,\mu_2$ and $\Omega_1,\Omega_2$, the performance can be nearly matched with the oracle if the posterior distributions of $\mu_1,\mu_2$ and $\Omega_1,\Omega_2$ contract near the true values sufficiently fast. In the high-dimensional setting, this is possible only if there is some lower-dimensional structure, such as sparsity. Du and Ghosal \cite{du} considered a prior on a precision matrix $\Omega$ based on a sparse modified Cholesky decomposition $\Omega = LDL^{\mathrm{T}}$, where $L$ is sparse and lower triangular with diagonal entries $1$ and $D$ is a diagonal matrix. The greatest benefit of using a Cholesky decomposition in the sparse setting is that the positive definiteness restriction on $\Omega$ is automatically maintained. However, a drawback is that the sparsity levels in the rows of $\Omega$ depend on the coordinate ordering, which may be somewhat arbitrary. To alleviate the problem, with a hard-spike-and-slab prior for the entries of $L$, the probability of a non-zero entry should be decreased with the row index $i$ proportional to $i^{-1/2}$. Then it follows that, when $i$ and $j$ are roughly proportional and both large, the probabilities of a non-zero in the $i$th and $j$th rows of $\Omega$ are roughly equal. Du and Ghosal \cite{du} used a soft-spike-and-slab or the horseshoe priors, for which also a stationary (approximate) sparsity in the off-diagonal elements of $\Omega$ can be maintained by making the slab probability decay like $i^{-1/2}$ for the former and the local parameter decay at this rate for the horseshoe prior. Du and Ghosal \cite{du} showed that the misclassification rate of the Bayes procedure converges to that of the oracle Bayes classifier for a general class of shrinkage priors when $p^2 (\log p)/n\to 0$, provided that the number of off-diagonal entries in the true $\Omega$ is $O(p)$. This is a substantial improvement over the requirement $p^4/n\to 0$ needed to justify this convergence without assuming sparsity of $\Omega$.  

\subsection{Exponential family with a graphical structure}

Representing various graphical models as exponential families facilitates inference for graph selection and structure learning. Suppose $X = (X_1,\ldots,X_p)$ is a $p$-dimensional random vector and $G=(V,E)$ is an undirected graph. The {\it exponential family graphical model} (with respect to $G$) has the joint density (with respect to a fixed dominating measure $\mu$ like the Lebesgue or the counting measure)  of the form
\begin{equation}
	p(X;\theta) = \exp\{\sum_{r \in V}\theta_rB(X_r) + \sum_{(r,t) \in E}\theta_{rt}B(X_r)B(X_t) + \sum_{r \in V}C(X_r) - A(\theta) \},
\end{equation}
for sufficient statistics $B(\cdot)$, base measure $C(\cdot)$ and log-normalization constant $A(\theta)$. Note that for the choice of $B(X) = X/\sigma, C(X) = -X^2/2\sigma^2$, we get the Gaussian graphical model
\begin{equation}
	p(X;\theta) \propto \exp\{\sum_{r \in V}\frac{1}{\sigma_r}\theta_rX_r + \sum_{(r,t) \in E}\frac{1}{\sigma_r \sigma_t}\theta_{rt}X_rX_t  - \sum_{r \in V}\frac{1}{\sigma_r^2}X_r^2 \}; 
\end{equation}
here $\{\theta_{rt}\}$ are the elements of the corresponding precision matrix.

An alternative formulation using a matrix form, called an {\it exponential trace-class model}, was recently proposed by Zhuang et al. \cite{zhuang2016graphical}. The family of densities for $X$  is indexed by a $q\times q$ matrix $M$, and is given by the expression 
$f(X|M)=\exp[-\langle M, T(X)\rangle +\xi(X)-\gamma(M)]$, 
where $\langle \cdot,\cdot\rangle$ 
stands for the trace inner product of matrices, and $T:\mathbb{R}^p \to \mathbb{R}^{q\times q} $, $q$ may or may not be the same as $p$, and $\gamma(M)$ is the normalizer for the corresponding exponential family given by 
$\gamma(M)=\log \int \exp [-\langle M, T(X)\rangle +\xi(X)]d\mu(X)$. The Gaussian graphical model is included as $q=p$, $T(X)=X X^{\mathrm{T}}$. In this case, $M$ agrees with the precision matrix and is positive definite. 
The most interesting feature in an exponential trace class model is that the conditional independence of $X_i$ and $X_j$ given others is precisely characterized by $M_{ij}=0$. In the high-dimensional setting, it is sensible to impose sparsity of the off-diagonal entries. In the Bayesian setting, this can be addressed by spike-and-slab or continuous shrinkage-type priors.

\subsubsection{Ising model}
For Bernoulli random variables defined over the nodes of the graph $G$, the exponential family representation of the graphical model, with the choice of $B(X) = X$ and the counting measure as the base measure, gives the distribution
\begin{equation}
	p(X;\theta) = \exp \{\sum_{r \in V}\theta_rX_r + \sum_{(r,t) \in E}\theta_{rt}X_rX_t  - A(\theta)\},
\end{equation}
which is popularly known as the Ising model. The conditional distribution of the nodes gives a logistic regression model. The model may also be represented as an exponential trace class model. Graph selection can be achieved via neighborhood-based variable selection method (Meinshausen and B\"uhlmann \cite{meinshausen2006high}) using the $\ell_1$-penalty. In the high-dimensional setting, Barber and Drton \cite{barber2015high} used the Bayesian Information Criterion (BIC) in the logistic neighborhood selection approach. A sparsity-based conditional MLE approach was proposed by Yang et al. \cite{yang2015graphical}. In the Bayesian setting, continuous shrinkage priors or spike-and-slab priors may be used to infer the graphical structure. Variational methods for inference on the parameters and evaluation of the log-partition function $A(\theta)$ were discussed in Wainwright and Jordan \cite{wainwright2008graphical}.

\subsubsection{Other exponential family graphical models}

Count data are often obtained from modern next-generation genomic studies. To study relations in count data, variables can be modeled via a Poisson distribution on the nodes, giving the Poisson graphical model $$p(X;\theta) = \exp \{\sum_{r \in V}(\theta_rX_r - \log(X_r!)) + \sum_{(r,t) \in E}\theta_{rt}X_rX_t  - A(\theta)\}.$$ Integrability forces $\theta_{rt} \leq 0$, so that only negative conditional relationships between variables are captured. 
An alternative formulation of a Poisson graphical model as an exponential trace-class model can accommodate positive interactions between variables, where the marginal distributions are Poisson but the conditional distributions are not (Zhuang et al. \cite{zhuang2016graphical}). A Poisson graphical model in which interactions are modeled via a product of powers of arctan transformations of the observations was studied by Roy and Dunson \cite{roy2020nonparametric}. 
Other useful exponential families include the multinomial Ising model.

\subsection{Nonparanormal graphical model} 

A semiparametric extension of GGMs that can model continuous variables located at nodes of a graph is given by the nonparanormal model (Liu et al. \cite{nonparanormal}), where it is assumed that the vector of variables $X = (X_1,\ldots,X_p)\in [0,1]^p$ reduces to a multivariate normal vector through $p$ monotone increasing transformations: for some monotone functions $f_1,\ldots,f_p$, 
$$f(X) := (f_1(X_1),\ldots, f_p(X_p))\sim \mathrm{N}_p(\mu, \Sigma)$$ 
for some $\mu\in \mathbb{R}^p$ and positive definite matrix $\Sigma$. The model is not identifiable and needs to fix the location and scale of the functions or the distribution. 
Liu et al. \cite{nonparanormal} developed a two-step 
estimation process in which the functions were estimated first using a truncated empirical distribution through the relations $f_j(x) =\Phi^{-1} (F_j (x))$, where $F_j$ stands for the cumulative distribution function of $X_j$. 
Mulgrave and Ghosal \cite{mulgrave1,mulgrave2,mulgrave3} considered two different Bayesian approaches --- based on imposing a prior on the underlying monotone transforms in the first two papers, and based on a rank-likelihood which eliminates the role of the transformations in the third. In the first approach, a finite random series based on a B-spline basis expansion is used to construct a prior on the transformation. The advantage of using a B-spline basis is that, to maintain monotonicity, the coefficients only need to be increasing. A multivariate normal prior truncated to the cone of ordered values can be conveniently used. However, to ensure identifiability, two constraints  
$$f(0)=1/2 \quad \mbox{ and } \quad f(3/4)-f(1/4)=1$$ 
are imposed, which translate to linear constraints, and hence the prior remained multivariate normal before imposing the order restriction. Samples from the posterior distribution of the ordered multivariate normal coefficients can be efficiently obtained using the exact Hamiltonian MCMC (Packman and Paninski \cite{HMCMC}). In Mulgrave and Ghosal \cite{mulgrave2}, a normal soft-spike-and-slab prior was put on the off-diagonal elements of $\Omega=\Sigma^{-1}$, and the scaling-it-up technique (Wang \cite{wang2015scaling}) was utilized to avoid the evaluation of the intractable normalizing constant arising from the positive definiteness restriction. They also proved a consistency result for the underlying transformation in terms of a pseudo-metric given by the uniform distance on a compact subinterval of $(0,1)$. The approach in Mulgrave and Ghosal \cite{mulgrave2} used a connection to the problem of regressing a component on subsequent components, where the partial regression coefficients form the basis of a Cholesky decomposition of the precision matrix. Sparsity in these coefficients was introduced through continuous shrinkage priors by increasing sparsity with the index as in Subsection~\ref{discriminant}. The resulting procedure is considerably faster than the direct approach of Mulgrave and Ghosal \cite{mulgrave1}. Mulgrave and Ghosal \cite{mulgrave2} also considered the mean-field variational approach for even faster computation. In Mulgrave and Ghosal \cite{mulgrave3}, the posterior distribution is altered by replacing conditioning on the data by conditioning on the ranks, which are the maximal invariants under all monotone transformations. The benefit is that the likelihood is free of transformations, eliminating the need to assign prior distributions to them and allowing them to be treated as fixed, even though they are unknown. The absence of the nonparametric part of the parameter thus makes the procedure a lot more efficient and robust. The posterior can be computed via Gibbs sampling using a simple data augmentation technique applied to the transformed variables. Also, the arbitrary centering and scaling means that only a scale-invariant of $\Omega$ is identifiable. Mulgrave and Ghosal \cite{mulgrave3} also derived a consistency result in a fixed-dimensional setting using Doob's posterior consistency theorem (see Ghosal and van der Vaart \cite{gvbook17}). 

An alternative way of expressing a joint distribution is by using the notion of a copula function $C(x_1,\ldots,x_p)$, which is a joint distribution function with uniform marginals and unspecified one-dimensional distribution functions $F_1,\ldots,F_p$ on $\mathbb{R}$  (Trivedi and Zimmer \cite{trivedi2007copula}).   
It may be noted that the structure of nonparanormality may be equivalently expressed by a  Gaussian copula model with 
$$C(x_1,\ldots,x_p)=\Phi_p (\Phi^{-1}(x_1),\ldots,\Phi^{-1}(x_p);R),$$ 
where $\Phi_p(\cdot;R)$ is the joint multivariate normal distribution with mean zero and a given covariance matrix $R$, by taking $C$ the joint distribution function of $(\Phi(f_1(X_1)-\mu_1),\ldots,  \Phi(f_p(X_p)-\mu_p))$. Pitt et al. \cite{pitt2006efficient} constructed a prior distribution and proposed an efficient posterior sampling technique in a Gaussian copula setting to make inference about conditional dependence.  A convenient method for creating a multidimensional copula is to use a vine copula, which combines pairwise distributions to form a joint distribution, allowing for the introduction of conditional independence sequentially as the prior is constructed; see Gruber and Czado \cite{gruber2015sequential} for a method and the corresponding posterior-computing technique. 

The nonparanormal model has also been used for (possibly high-dimensional) discriminant analysis by Zhu et al. \cite{zhu2025bayesian} in the semisupervised learning setting, where the population labels of most observations are missing. 

\section{Estimating a long vector smoothly varying over a graph}
\label{sec:smooth graph}

Let us revisit the normal sequence model $X_i=\theta_i+\varepsilon_i$, $i=1,\ldots,n$, where $\varepsilon_i \sim \mathrm{N}(0,\sigma^2)$ independently, as in Section~2, but the sparsity assumption on the vector of means $f=(\theta_1,\ldots,\theta_n)$ is not appropriate. Instead, the values are assumed to `smoothly vary over locations' in some appropriate sense. The simplest case is that these values lie on a linear lattice; more generally, the positions may correspond to the nodes of a graph. Smooth variation of $\theta_i$ with respect to $i$ over a graph can be mathematically quantified through the graph Laplacian $L=D-A$, where $D$ is the diagonal matrix of node degrees, and $A$ stands for the adjacency matrix. We define H\"older ball of $\beta$-smoothness of radius $Q$ over a graph of size $n$ with graph Laplacian $L$ by 
$$ \mathcal{H}_\beta(Q)=\{f=(\theta_1,\ldots,\theta_n):\langle f, (I+(n^{2/r}L)^\beta) f\rangle \le Q^2\},$$ 
where $r$ stands for the ``dimension'' of the graph defined through the growth of eigenvalues of the graph Laplacian $L$: the $i$th eigenvalue grows like $(i/n)^{2/r}$. For lattice graphs, the dimension truly agrees with the physical dimension, but in general, it can be a fractional number. In this setting, Kirichenko and van Zanten \cite{KirichenkoMinimax} showed that the minimax rate of recovery on $\mathcal{H}_\beta(Q)$ is $n^{-\beta/(2\beta+r)}$. Kirichenko and van Zanten \cite{Kirichenko} developed a Bayesian procedure for this problem using a multivariate normal prior: 
$$f\sim \mathrm{N}_n (0, (n/c)^{(2\alpha+r)/r} (L+ n^{-2}I)^{(2\alpha+r)/r}),$$ 
where $c$ is given the standard exponential prior. Using van der Vaart-van Zanten posterior contraction rate theory for a  Gaussian process prior (Ghosal and van der Vaart \cite{gvbook17}, Chapter 11), they showed that the posterior contracts at the optimal rate  $n^{-\beta/(2\beta+r)}$ up to a logarithmic factor, whenever the true smoothness is $\beta\le \alpha+r/2$. An analogous result also holds for binary regression. Kirichenko and van Zanten \cite{Kirichenko} also showed that the limited range adaptation can be strengthened to full range adaptation using an exponential of the Laplacian covariance kernel. 

An extension of the setup of Kirichenko and van Zanten \cite{Kirichenko} for functional data, with the means taking values in a Hilbert space, was considered by Roy and Ghosal \cite{ArkaGraph}. They introduced a notion of Sobolev smoothness classes indexed by two smoothness indices: $\beta$ for graphical smoothness and $\gamma$ for functional smoothness. For a model with the graphical dimension $r$, their obtained rate $n^{-\beta\gamma/(2\beta\gamma+\beta+r\gamma)}$ is an upper bound for the minimax rate of recovery in the mixed Sobolev smoothness class and is minimax over the slightly larger homogeneous anisotropic smoothness class; see Belitser and Ghosal \cite{belitser2025mean} for the definitions of these smoothness classes and a clarification of Roy and Ghosal's \cite{ArkaGraph} minimaxity result as an application of their treatment of multidimensional infinite arrays under various smoothness settings. Roy and Ghosal \cite{ArkaGraph} also designed Bayes procedures that achieve this rate exactly, and also adaptively within a logarithmic factor. In a non-adaptive setting, they further showed that a slightly inflated Bayesian credible ball of optimal size has adequate frequentist coverage.

\section{Matrix models}

Several high-dimensional problems involve unknown matrices instead of vectors. In Section~5, we discussed the literature on Bayesian inference for the covariance and precision matrices of high-dimensional observations. In this section, we review results when observations are large matrices. Some examples include multiple linear regression with group sparsity, multi-task learning, matrix completion, stochastic block model, and biclustering.

\subsection{Generic results in structured linear models}

Gao et al. \cite{gvz15} derived results on posterior dimensionality and contraction rates for many models simultaneously in a general `signal plus noise' model when the signal is a vector or a matrix having some specific structure. One example is multiple linear regression with group sparsity, where one observes $Y=XB+W$, with $X$ being an $n\times p$ matrix, $B$ a $p\times m$ matrix, and the columns of $B$ are assumed to share common support of size $s$. They obtained the minimax rate $s(m+\log(ep/s))$ for the prediction problem of estimating $XB$ in the Frobenius norm. This model is a special case of multi-task learning in which the columns of $B$ share a specific structure. For instance, instead of assuming a joint sparsity pattern among columns, one may think that columns of $B$ can only be chosen from a given list of $k$ possible columns, with $k$ typically much smaller than $m$. The posterior rate $pk+m\log{k}$ was derived for prediction in Gao et al. \cite{gvz15} (it is optimal as soon as  $k<pm$). Dictionary learning, on the other hand, assumes that the signal matrix of size $n\times d$ is $\theta=QZ$, for $Q$ an $n\times p$ dictionary matrix and $Z$ a discrete matrix of size $p\times d$ with sparse columns. Gao et al. \cite{gvz15} derived the adaptive rate $np+ds\log(ep/s)$ for estimating $\theta$, which is optimal if smaller than $nd$. Kim and Gao \cite{kimgao19} proposed an EM-type algorithm to simulate from posterior aspects corresponding to the prior considered in Gao et al. \cite{gvz15}, recovering as a special case the EMVS algorithm of Ro\u ckov\`a and George \cite{rockovageorge14} for linear regression. Belitser and Nurushev \cite{belitsernur19} also followed this general framework, considering a re-centered normal prior extending the approach in Belitser and Nurushev \cite{belnur15} and Belitser and Ghosal \cite{belitserghosal16}, and derived both local oracle posterior rates, as well as optimal-size credible balls with guaranteed coverage.

\subsection{Stochastic block model and community detection}

In a \emph{stochastic block model} (SBM) with $k$ groups, one observes  $Y = \theta + W$ with $\theta_{ij}=Q_{z(i)z(j)}\in [0,1]^{n\times n}$ for some matrix $Q$ of size $k\times k$ of edge probabilities, a labeling map $z\in\{1,\ldots,k\}^n$ and $W$ a centered Bernoulli noise. Gao et al. \cite{gvz15} treated the question of the estimation of $\theta$ in the Frobenius norm within their general framework, getting the adaptive minimax rate $k^2+n\log{k}$ with unknown $k$. A similar result and coverage of credible sets were derived by Belitser and Nurushev \cite{belitsernur19}.  Pati and Bhattacharya \cite{patibhatt15} obtained the near-optimal rate $k^2\log ({n/k})+n\log{k}$ for known $k$ using an independent uniform prior on the coordinates of $Q$ and a Dirichlet prior probabilities for the label proportions. The biclustering model can be viewed as an asymmetric extension of the SBM model, where $\theta$ is an $n\times m$ rectangular matrix, and the rows and columns of $Q$ have their own labels, with $k$ and $l$ groups, respectively. This model was handled in Gao et al. \cite{gvz15} and Belitser and Nurushev \cite{belitsernur19} using similar techniques as before, with the optimal adaptive posterior rate $kl+n\log{k}+m\log{l}$.

Another popular question is that of the recovery of the vector of labels $z$ (possibly up to a permutation), also known as \emph{community detection}. This is a more ambitious task than estimating $\theta$. It can only be guaranteed either asymptotically (Bickel et al. \cite{bickel13}), or non-asymptotically, imposing some separation between classes (Lei and Zhu \cite{Lei:Zhu:2017:1}); see Castillo and Orbanz \cite{co17} for a discussion on uniformity issues. Van der Pas and van der Vaart \cite{vpvv18} showed that the posterior mode corresponding to a beta prior on edge probabilities and Dirichlet prior probabilities for the label proportions asymptotically recovers the labels when the number of groups $k$ is known, provided that the mean degree of the graph is at least of order $\log^2 n$. Their approach relies on studying the \emph{Bayesian modularity}, that is, the marginal likelihood of the class labels given the data, with the edge probabilities integrated out. Kleijn and van Waaij \cite{kleijn2018recovery} and van Waaij and Kleijn \cite{waaij2020uncertainty} derived estimation and uncertainty quantification results for posterior distributions in the special case of the planted multi-section model, where the matrix $Q$ has only two different values, one for diagonal elements and another for off-diagonal ones.

\subsection{Matrix completion}

In the matrix completion model (also known as the ``Netflix problem''), one observes $n$ noisy entries of an unknown $m\times p$ matrix $M$, typically assumed to be of a low rank $r$ (or well-approximated by a low-rank matrix).

% The entries are sampled at random from a distribution on entry indices. The rate $(m+p)r\log \max(m, p)$, minimax optimal up to the logarithmic term, was obtained by Mai and Alquier \cite{maialquier15} in the PAC-Bayesian setting for a prior sitting close to small-rank matrices and a class of tempered posteriors (see the discussion section below for more on this topic) while Suzuki \cite{suzuki15} derived similar results for posterior distributions (under assumptions on the law of the design) and also obtained results in the more general case of tensors. [Please add more descriptions of the procedures]

A convenient formulation is
\[
Y_i=M_{\omega_i}+\varepsilon_i,\qquad i=1,\ldots,n,
\]
where $\omega_i=(j_i,k_i)\in\{1,\ldots,m\}\times\{1,\ldots,p\}$ is the observed entry index and $\varepsilon_i$ is an observation error. The indices $\omega_i$ are usually sampled from a distribution on the set of matrix entries. The uniform design is the simplest case, but in many applications, including recommender systems, the sampling distribution is far from uniform. The inferential objective is to recover $M$, or to predict its unobserved entries, under the structural assumption that $M$ has low rank, or is well approximated by a low-rank matrix.

A Bayesian formulation typically introduces the low-rank structure through a factor representation \(M=UV^{\mathrm{T}}\), where $U\in\mathbb R^{m\times K}$ and $V\in\mathbb R^{p\times K}$ for a working rank $K$. The working rank may be chosen as an upper bound, while the prior is designed to shrink unnecessary factors. For example, one may place hierarchical Gaussian priors on the columns of $U$ and $V$, with column-specific scale parameters that determine the effective rank. If the scale corresponding to a given column is strongly shrunk toward zero, the associated rank-one component
$U_{\cdot \ell}V_{\cdot \ell}^{\mathrm{T}}$ contributes little to the reconstructed matrix. This is the Bayesian analog of rank regularization. Probabilistic matrix factorization and Bayesian probabilistic matrix factorization apply this idea directly to collaborative filtering, treating the rows and columns as latent user and item features and learning their posterior distributions from observed ratings (Salakhutdinov and Mnih \cite{salakhutdinov08}). More general sparse Bayesian low-rank formulations also use hierarchical shrinkage priors to control the effective rank and to regularize the factor matrices (Babacan et al. \cite{babacan12}).

The likelihood is usually specified only on the observed entries. In the Gaussian case, this gives 
\[
p(Y\mid U,V,\sigma^2,\omega)
=
\prod_{i=1}^n
\phi_{\sigma}
(Y_i-U_{j_i}^{\top}V_{k_i}),
\]
where $U_{j_i}$ and $V_{k_i}$ denote the corresponding row vectors of the two factor matrices, and $\phi_\sigma$ is the centered normal density with standard deviation $\sigma$. The posterior then combines this observed-entry likelihood with priors on $U$, $V$, $\sigma$, and possibly the effective rank. Prediction of an unobserved entry $(j,k)$ is obtained through the posterior
predictive distribution of $U_j^{\mathrm{T}} V_k$, or through the posterior mean
\[
\widehat M_{jk}
=
\int U_j^{\mathrm{T}} V_k \, d\Pi(U,V\mid Y,\omega).
\]
Thus, unlike a point estimator based on the nuclear-norm penalization, the Bayesian procedure gives a full distribution over possible completed matrices and therefore quantifies uncertainty in the missing entries.

Mai and Alquier \cite{maialquier15} studied this problem from a PAC-Bayesian perspective. Their procedure is based on a Gibbs, or tempered, posterior of the form
\[
\Pi_\lambda(dA\mid Y,\omega)
\propto
\exp\big\{
-\lambda \sum_{i=1}^n (Y_i-A_{\omega_i})^2
\big\}
\Pi(dA),
\]
where $\lambda>0$ is an inverse-temperature parameter and $\Pi$ is a prior concentrating near low-rank matrices. This is a generalized posterior obtained by exponentiating the empirical
loss, and not necessarily the ordinary Bayesian posterior corresponding to a correctly specified Gaussian likelihood. The posterior mean under this Gibbs distribution is then used as an aggregate estimator. Their analysis shows that, with a prior suitably adapted to low-rank structure, the resulting procedure satisfies oracle inequalities and attains the rate \((m+p)r\log\{\max(m,p)\},\) up to constants and logarithmic factors. This rate is minimax optimal up to the logarithmic term. The PAC-Bayesian argument is well-suited to this setting because the complexity of the estimator is controlled through the Kullback--Leibler divergence between a localized posterior distribution and the prior, rather than through an explicit combinatorial count of low-rank models.

Suzuki \cite{suzuki2015convergence} developed related Bayesian concentration results for matrix (and also, tensor) completion. In the matrix case, the posterior is again built from an observed-entry likelihood and a prior favoring low-rank structure, but the analysis is expressed as posterior concentration around the true matrix. A key feature of this work is that the sampling distribution of the entries is explicitly reflected in the loss under which contraction is measured. Under appropriate assumptions on the design distribution, the posterior contracts at a rate governed by the effective number of free parameters in a low-rank matrix. 

From the Bayesian point of view, matrix completion is therefore a useful prototype for posterior contraction in models where the ambient dimension $mp$ is much larger than the intrinsic dimension. The prior encodes the low-rank geometry, while the likelihood uses only the randomly observed entries. Tempered or generalized posteriors, as in the PAC-Bayesian analysis of Mai and Alquier \cite{maialquier15}, provide one route to sharp oracle inequalities;
ordinary or hierarchical Bayesian posteriors, as in  Suzuki \cite{suzuki2015convergence}, provide a direct concentration theory. Both approaches show that, under suitable low-rank priors and regularity assumptions on the sampling design, Bayesian procedures can achieve statistical rates essentially the same as those of the best-known frequentist methods, while retaining a posterior distribution over the completed matrix.

\subsection{Matrix regression}

The same low-rank Bayesian principle also appears in matrix regression. In this setting, one observes
\[
Y_i=\langle X_i,M\rangle+\varepsilon_i,\qquad i=1,\ldots,n,
\]
where $X_i\in\mathbb R^{m\times p}$ is a matrix-valued covariate and $M$ is an unknown coefficient matrix. Matrix completion is obtained as the special case $X_i=e_{j_i}e_{k_i}^{\mathrm{T}}$, so that $\langle X_i,M\rangle=M_{j_i k_i}$. More general matrix regression allows the observation operator to be dense or structured, and the statistical difficulty is again controlled by the effective dimension of a low-rank matrix rather than by the ambient dimension $mp$. A Bayesian treatment places a prior on $M$ through a factorization $M=UV^{\mathrm{T}}$, or through a prior on its singular values and singular vectors, and combines it with the regression likelihood. Posterior contraction is then expected to depend on the intrinsic dimension $(m+p)r$, up to logarithmic factors and constants depending on the design. This connects matrix completion with the broader class of trace-regression or matrix-sensing problems studied in the low-rank estimation literature (Negahban and Wainwright \cite{negahban11}, Rohde and Tsybakov \cite{rohde11}).

\section{Tensor models}

 Tensors are higher-dimensional arrays that generalize matrices and arise in many modern applications, such as context-aware recommender systems, neuroimaging, genomics, spatio–temporal modeling, and multi-relation networks. We provide a brief overview of Bayesian methods addressing these problems.

\subsection{Tensor completion problem}

Suppose that we observe certain entries, to be collectively called $\Omega$, of an order $d$-tensor $\mathcal{Y}$ which has mean  $\Theta\in\mathbb{R}^{n_1\times\ldots\times n_d}$. The problem is to estimate the mean tensor $\Theta$ and thus predict the unobserved entries. When $\Omega$ is the whole index set, the problem simply reduces to denoising. Belitser and Ghosal \cite{belitser2026bayesian} considered the denoising problem for infinite-length tensors, assuming a variety of smoothness conditions given by the decay rates of the entries, and obtained oracle rates of concentration relative to a family of structures for an empirical Bayes posterior distribution assuming a normal working model for the noise and normal priors on the coefficients. They allowed different smoothness parameters in different directions, as well as their dependence on a leading index, and showed that, for recovery at the oracle rate in general, one must consider the `inhomogeneous family of structures' consisting of collections of index-vectors $(i,j_1,\ldots,j_d)$ of the form $i\le I$, $j_k\le J_{k,i}$ depending on both $i$ and $k$, $k=1,\ldots,d$. They illustrated their rates with several examples, including the models discussed in Section~\ref{sec:smooth graph}, and recovered the stated rates.

As in matrix completion, a low-rank structure is often assumed. However, there are multiple notions of tensor ranks. The popularly used Tucker decomposition  (Tucker \cite{tucker1966some}) of a tensor is given by 
\begin{align*} 
\Theta &= \mathcal{C} \times_1 U_1 \times_2 U_2 \times_3 \ldots \times_d U_d\\
&:=(\sum_{j_1}^{m_1}\cdots 
\sum_{j_d=1}^{m_d} c_{j_1\ldots j_d} u_{i_1 j_1}\cdots u_{i_d j_d}: i_k=1,\ldots,n_k,\, k=1,\ldots,d)
\end{align*} 
where $\mathcal{C} \in \mathbb{R}^{m_1\times \ldots \times m_d}$ is called the core tensor and $U_k \in \mathbb{R}^{n_k \times m_k}$ is the $k$th factor matrix, $k=1,\ldots,d$; here . The Tucker (multilinear) rank of $\mathcal{Y}$ is said to be $(m_1, \ldots, m_d)$. 
The assumption reduces the number of parameters from $\prod_{k=1}^d {n_k}$ to a much smaller number $\prod_{k=1}^d {m_k}+\sum_{k=1}^d n_k m_k$, since $m_k\ll n_k$ for all $k=1,\ldots,d$. Additional sparsity conditions on the factor matrices $U_1,\ldots, U_d$ may be imposed. Only recently has the Bayesian tensor completion problem been addressed. Chu and Ghahramani \cite{chu2009probabilistic} extended the technique of Bayesian matrix factorization to tensors, and their method can automatically infer the tensor ranks. Xiong et al. \cite{xiong2010temporal} introduced the Bayesian Probabilistic Tensor Factorization (BPTF) method for recommender systems by imposing a prior distribution on user-item-context interactions. Further developments took place 
through the works of Liu et al. \cite{liu2012tensor}, Zhao et al. \cite{zhao2015bayesian}, and Suzuki \cite{suzuki2015convergence}. 
Suzuki \cite{suzuki2015convergence} showed that posterior concentration can be obtained in such tensor models, with rates governed by the number of free parameters in the chosen low-rank tensor representation rather than by the full ambient tensor dimension. 
Rai et al. \cite{rai2014scalable} imposed a structured sparsity condition to obtain interpretable latent structures. Rai et al. \cite{rai2015leveraging} improved recovery under additional structural assumptions. Ghosh and Ghosal \cite{shuvrarghya_completion} imposed further structure on the factor matrices through smooth functions and latent variables of the form $f_{j_k}(\eta_{i_k})$ for the $(i_k,j_k)$th entry of $U_k$. They considered priors based on random B-spline series and obtained the posterior contraction rate $(\prod_{k=1}^d m_k+\sum_{k=1}^d m_k n_k^{1/(1+2\alpha_k)})^{1/2}(N/\log N)^{-1/2}$ with $N=\prod_{k=1}^d n_k$, when the functions in the $k$th latent factor are H\"older $\alpha_k$-smooth. 
 Thus, matrix and tensor completion are special cases of a broader Bayesian theory for low-rank structured regression under partial or indirect observation.

\subsection{Tensor-on-tensor regression}

Suppose that we observe $N$ independent pairs of tensors $\{\mathcal{X}_i, \mathcal{Y}_i: i=1,\ldots,N\}$, respectively of order 
${n_1\times\ldots\times n_p}$ and 
${m_1\times\ldots\times m_q}$, related by the tensor-on-tensor regression model 
$$\mathcal{Y}_i = \langle \mathcal{X}_i, \mathcal{B}  \rangle_p + \mathcal{E}_i,$$ 
with error $\mathcal{E}_i$, $i = \{1,2,\ldots,N\}$, 
where $\langle \mathcal{C}, \mathcal{D} \rangle_p$ denotes the $p$-fold the \textit{tensor contracted product} of 
$$\mathcal{C}=(c_{i_1\cdots i_p}: i_k=1,\ldots, m_k,\, k=1,\ldots,d\}$$ and 
$$\mathcal{D} =(d_{i_1\cdots i_p,j_1\cdots j_q}:i_k=1,\ldots, n_k,\, k=1,\ldots,p,\;  j_l=1,\ldots, m_l,\, l=1,\ldots,q)$$ given by $$(\sum_{i_1}^{n_1}\cdots \sum_{i_p=n_p} c_{i_1\cdots i_p} d_{i_1\cdots i_p,j_1\cdots j_q} :j_l=1,\ldots, m_l,\, l=1,\ldots,q),$$ and the coefficient tensor $\mathcal{B} \in \mathbb{R}^{n_1\times\ldots\times n_p \times m_1 \times\ldots\times m_q}$. The ambient dimension of $\mathcal{B}$ is huge: $\prod_{k=1}^p n_k\prod_{k=1}^q m_k$, typically much higher than $N$. For consistent estimation, structural assumptions are imposed on $\mathcal{B}_0$. 
The common special case of scalar-on-tensor regression 
is obtained when $\mathcal{Y}_i$'s are scalar, while tensor-on-scalar regression corresponds to scalar $\mathcal{X}_i$'s. 

Lock \cite{Lock2018} is apparently the first to study the tensor-on-tensor model within the Bayesian paradigm. Lock \cite{Lock2018} imposed the condition that the coefficient tensor has a CP decomposition, which is a sum of rank-one tensor factorizations and is a special case of the Tucker decomposition, and to estimate the latent factors, used a ridge regularization technique. Wang and Xu \cite{wang2024} approached the  Bayesian regression problem based on a Tucker decomposition of the coefficient tensor. Under mild conditions, a posterior consistency result for the scalar-on-tensor regime was obtained by Guhaniyogi et al. \citep{guhaniyogi2017}, and Papadogeorgou et al. \cite{papadogeorgou2021}. The posterior convergence theory for the general tensor-on-tensor models was unexplored until recently, when Ghosh and Ghosal \cite{shuvrarghya_regression} considered the general model with a Tucker decomposition for the coefficient tensor, with the core tensor $\mathcal{B} = \mathcal{C}\times_1 U_1\times_2\ldots\times_{p+q}  U_{p+q}$, 
having the multi-linear rank 
$(r_1, \ldots,r_{p+q})$,  and assumed a low-rank plus sparse model for each matrix factor $U_k$. They derived the posterior contraction rate  $(\prod_{k=1}^{p+q} r_k +\sum_{i=1}^p n_i r_i+\sum_{j=1}^q m_j r_{p+j})^{1/2} (N/\log N)^{-1/2}$. They also obtained a Bernstein-von Mises theorem for the coefficient in a vectorized form.

\section{High-dimensional time series}

Thus far, our discussions have mostly concerned independent data. Observations over time show dependence, which must be addressed in the modeling. Data that measure various characteristics simultaneously are often high-dimensional and omnipresent in economic time series. Standard time series models, such as autoregressive (AR) or autoregressive moving average (ARMA) processes, allow straightforward multivariate generalizations, and desirable properties such as stationarity and causality can be incorporated in the model. In addition to making inferences about parameters that drive the time series, detecting the intrinsic dependence structure among time series components, as in graphical models, is of the utmost importance and often the primary interest, with the time series structure merely a tool for accommodating sample correlation. 
When the samples are correlated, dynamic models, such as a stationary GGM, extend static GGM models to account for correlation via a multivariate time-series structure. 
Dynamic Factor Models (DFM), widely used in economic applications, are well-suited for estimating latent, lower-dimensional processes that drive dynamics, but their parameterization is unsuitable for estimating graphical and stationary structures; see Matteson and Tsay \cite{matteson2011dynamic} for a recent survey. 
Time series on graphs and networks were studied in Basu et al.  \cite{basu2015network} and Ma and Michailidis \cite{ma2016joint} from a frequentist consideration. Kolar et al. \cite{kolar2010estimating} and Chen et al. \cite{chen2013covariance} estimated graphical dependence using a sparse precision matrix, but their approaches cannot preserve stationarity and causality in the resulting processes.

Since the parameters governing the graphical structure and the temporal dynamics are intertwined, the likelihood and subsequent posterior sampling via MCMC in a Bayesian approach are considerably more complicated than the independent case. To alleviate the problem, Roy et al. \cite{roy2024bayesian} 
proposed a semiparametric model closely resembling the DFM, by representing the multivariate time-series $X(t)=(X_1(t),\ldots,X_p(t))^{\mathrm{T}}$ in terms of a linear combination of an orthogonal rotation $U$ of a vector of latent independent univariate stationary time series $Z(t)=(Z_1(t),\ldots,Z_p(t))^{\mathrm{T}}$ with unspecified smooth spectral densities: $X(t)=\Sigma^{1/2} U Z(t)$. The parameterization of the contemporaneous precision matrix $\Omega=\Sigma^{-1}$ allows graph estimation via the stationary precision matrix. The framework is reminiscent of Independent Component Analysis (ICA), but Roy et al. \cite{roy2024bayesian} primarily aimed to learn the graphical structure rather than to identify the latent source signals. They adopted a Bayesian approach and placed prior distributions by reparameterizing the precision matrix via a sparse Cholesky decomposition, the orthogonal rotation via the Cayley representation using a skew-symmetric matrix, and expanding the spectral densities in terms of B-spline series. Their constraint-free parameterization allows free Metropolis-Hastings movements for posterior sampling based on the Whittle likelihood. Assuming a lower-dimensional representation in the basis expansion of the spectral densities, they derived the posterior contraction rate for the precision matrix by slightly modifying the general posterior contraction rate theory (Ghosal and van der Vaart \cite{ghosal2017fundamentals}) to incorporate the misspecification of the likelihood by the Whittle approximation. In particular, they observed that, if the spectral densities are sufficiently smooth, then the posterior contraction rate for the precision matrix in the Frobenius norm agrees with the rate $\sqrt{((p+s_0)\log n)/n}$ for the i.i.d. case in Subsection~\ref{wo specific sparisty precision}, where $p$ is the dimension of the time series, $s_0$ is the number of non-zero off-diagonal elements of the true precision matrix, and $n$ is the number of regularly-spaced temporal observations. A slightly simpler alternative that replaces the latent processes with unspecified spectral densities by independent AR processes was proposed by Ghosh et al. \cite{ghosh2025bayesian}, and the posterior contraction rate of the i.i.d. data for the precision matrix was also obtained. 

Vector autoregressive (VAR) processes are widely used for multidimensional time series, but in high-dimensional settings, they have too many parameters to draw inference from limited-length time series unless a dimension-reduction strategy is employed, because strong temporal dependence greatly reduces the useful information content of a time series. Causality, which means the prediction formula includes only current and past innovations, stabilizing the forecast of the future based on the present, plays an important role in multivariate time series models. As causality imposes complex constraints on process parameters, a fruitful approach is to reparameterize in terms of unconstrained parameters to compute the likelihood and sample from the posterior distribution. For a $p$-dimensional VAR time series of order $m$, causality is equivalent to having all roots of a determinantal equation involving the VAR coefficients lying outside the unit disc. Roy et al. \cite{roy2019constrained} characterized the condition as 
${\Gamma}(0) = {C}_0  \geq C_1 \geq \cdots \geq {C}_m = {\Sigma}$, say, 
where  $\Gamma(0)$ is the $p\times p$-autocovariance matrix,  $C_j = \mathrm{Var}(X_{j+1} |X_{j}, \ldots, X_1)$, $j=0,\ldots,m$, are the conditional dispersion matrices. The condition is clearly equivalent to $\Omega \leq C_1^{-1} \leq \cdots \leq C_m^{-1} = \Sigma^{-1}$. Roy et al. \cite{roy2024relational} proposed a dimension reduction by expressing the differences as low-rank symmetric nonnegative definite matrices $L_j L_j^{\mathrm{T}}$, with $L_j$s being a constraint-free $p\times r_j$-matrix, $r_j\ll p$, $j=1,\ldots,m$. Specifically, rank-one updates ($r_1=\cdots=r_m=1$) allow very explicit formulas for recursive estimation of the log-likelihood function, and the resulting efficient posterior sampling based on a sparse Cholesky decomposition prior for $\Omega$. Roy et al. \cite{roy2024relational} showed that under some regularity conditions, the AR parameters are estimated at the rate $p\sqrt{(\log n)/n}$ in the Euclidean distance and the precision matrix $\Omega$ controlling the graphical structure of the component series is estimated at the rate $\sqrt{((p+s_0)\log n)/n}$ based on $n$ temporal observations. 

Understanding the underlying marginal causal relationships among variables in a multivariate time series may be of primary interest for some applications, especially in econometrics. A directed acyclic graph (DAG) can succinctly represent such causal dependencies. A Bayesian approach for a multivariate stationary time series with a DAG structure but not assuming any pre-specified parent-child ordering was proposed by Roy et al. \cite{roy2025bayesian}. With $W$ standing for the matrix decorrelating the components of the multivariate time series $X(t)$ in that the de-caused series $X(t) -W X(t)$ has independent components, Roy et al. \cite{roy2025bayesian} considered independent priors based on B-splines series expansion as in Roy et al. \cite{roy2024bayesian} on the spectral densities of those independent time series, and placed a sparse prior on $W$. The traditional posterior in this problem is extremely complicated, as all DAG structures must be considered a priori. Roy et al. \cite{roy2025bayesian} adopted the sparse projection-posterior method described in Subsection~\ref{sparse projection posterior}, and developed a posterior computational algorithm. The model has an obvious identifiability problem due to the possibility of label switching. Roy et al. \cite{roy2025bayesian} provided two convergence results in terms of the Frobenius norm: under some regularity conditions, (a) if the variances of the de-caused series are known (which ensures identifiability) and the spectral densities are sufficiently smooth, then the sparse projection posterior for $W$ concentrates at the rate $\sqrt{((p+s_0)\log n)/n}$, where $s_0$ is the number of non-zero entries of true $W$; (b) (in the absence of the knowledge of the variance) if the spectral densities are linear independent, the spare projection posterior concentrates in arbitrarily small fixed-sized neighborhoods of the true $W$. They also showed that their model extends to cover matrix-variate time series. 

\section{Further topics}

\subsection{Variational Bayes}

{Variational Bayes (VB) methods have gained much popularity recently as a possible way to compute certain posterior approximations efficiently. The VB method first chooses a class $\mathcal{Q}$ of computationally tractable distributions. For instance, $\mathcal{Q}$ can be chosen to be the set of distributions that make the coordinates independent, the so-called {\em mean-field} class; another choice of $\mathcal{Q}$ is the set of parametric distributions of a given dimension. Second,  one chooses a loss, very often the Kullback-Leibler divergence, leading to the random distribution 
	\[ \tilde\Pi_X = \underset{Q\in\mathcal{Q}}{\text{argmin}}\ \text{KL}(Q;\Pi(\cdot | X)).\]
	In other words, the variational posterior $\tilde\Pi_X$ is the best approximation from the class $\mathcal{Q}$ of the given posterior $\Pi(\cdot | X)$. For nonparametric models, a theory for convergence rates of variational posteriors and $\alpha$-posteriors ($\alpha<1$, corresponding to a likelihood in Bayes' formula raised to the power $\alpha$) was obtained in Alquier and Ridgway \cite{alquierridgway}, Zhang and Gao \cite{zhanggao20} (where a link between VB and a class of empirical Bayes procedures was made) and Yang et al. \cite{ypb20}. The work Ch\'erief-Abdellatief \cite{badrphd} investigated convergence of variational posteriors in several settings, including mixture models and deep neural nonparametric regression. 
 Wang and Blei \cite{wb19} showed that the VB-posterior distribution in fixed-dimensional parametric models satisfies the following version of the Bernstein-von Mises theorem: the mean-field VB-posterior is approximated by the Kullback-Leibler projection of the posterior onto normal models centered at the truth, and the VB-posterior mean is a consistent and asymptotically normal estimator. They also argued that their result holds beyond the mean-field VB setting. 

	In the context of high-dimensional linear regression of Section \ref{regression}, variational posteriors are very popular and were considered among others in Carbonetto and Stephens \cite{carbonettostephens12} and Ormerod et al. \cite{ormerodetal17}, taking a set of spike-and-slab distributions as a variational class. Ray and Szab\'o \cite{rayszabo19} have analyzed the convergence rate of the corresponding VB--posterior based on a prior with Laplace slabs with parameter $\lambda$ and the mean-field variational class
	\[ \mathcal{Q}=\left\{ \bigotimes_{i=1}^p \{\gamma_i \mathrm{N}(\mu_i,\sigma_i^2)+(1-\gamma_i)\delta_0\},\ \ (\mu_i,\sigma_i^2,\gamma_i)\in \mathbb{R}\times \mathbb{R}^+\times [0,1]\right\}.\] 
	One can obtain the same convergence rates and variable selection properties as for spike--and--slab posteriors, under similar conditions for the parameter $\lambda$ and the design matrix $X$. Laplace slabs should be preferred to Gaussian slabs in the prior distribution (but not so in $\mathcal{Q}$), for the same reasons as for plain spike-and-slab.
	The computation of the VB posterior is an optimization problem that can be solved using the coordinate ascent variational inference (CAVI) algorithm: the parameters $\gamma_i, \mu_i, \sigma_i^2$ are sequentially updated, with the posterior approximated at each step by the variational subclass that keeps all parameters but one fixed. The main advantage of VB is computational speed. The procedure's tuning needs some care. However, the slab parameter $\lambda$ needs to be chosen, and the output of the CAVI algorithm is susceptible to the order in which parameters are updated. Ray and Szab\'o \cite{rayszabo19} proposed a prioritized update scheme to handle this issue. Similar results for VB--posteriors for the logistic regression were derived by Ray et al. \cite{rsc20}. If the goal is inference on a finite subset of coordinates (such as a given coordinate, e.g., $\beta_1$, or any finite subset of $\beta_i$'s) Castillo et al.  \cite{clrt24} proposed a VB class inspired by the reparameterization considered by Yang \cite{yang19}, derive asymptotic normality of the marginal VB posterior for the coordinates of interest in the form of a Bernstein-von Mises theorem if the covariates are not too correlated and provide empirical evidence under both low and high correlation among covariates. 
	 
	%[include other settings/examples of variational Bayes?]
	%This is equivalent to the maximisation of the evidence lower bound
	%\[ \tilde\Pi_X = \underset{Q\in\mathcal{Q}}{\text{argmax}}\ \int \ell_n(\theta;X)\mu(d\theta) - KL(Q||\Pi),
	%\]
}

\subsection{Link between Bayes and PAC-Bayesian approaches}

A popular approach in machine learning is to use a pseudo- (or generalized) Bayes measure by replacing the likelihood with a small power (also called the `temperature') in the standard Bayes posterior distribution. This can attenuate the dependence in the data, providing more robustness in misspecified settings. The log-likelihood may also be replaced by an empirical risk, depending on the loss function of interest. PAC-Bayesian inequalities enable to derive concentration of such pseudo-posterior measures around the true (or pseudo-true) parameter, for instance through upper bounds on the loss integrated over the (pseudo-)posterior. Seminal contributions to PAC-Bayes theory include those by McAllester, Catoni, and T. Zhang. Theory and algorithms for high-dimensional data were considered among others by Audibert, Alquier, Dalalyan and Tsybakov,  Martin and Walker, and Gr\"unwald and co-authors. We refer to the surveys by Guedj \cite{guedj19} and Alquier \cite{alquierpac24} for recent overviews of the field. We note that although such methods are sometimes referred to as generalized Bayes, they often assume that the temperature parameter is strictly less than $1$ (or tends to $0$ in some cases), thereby separating the corresponding results from those for original Bayes posteriors. Developing a unified theory across different temperature ranges, including the original Bayes, is an interesting direction for future work.

\subsection{Sparse projection posterior approach}
\label{sparse projection posterior}

The primary reason the posterior distribution in a sparse setting is complicated is that it must encode the sparsity-generating mechanism in the prior — either via a spike-and-slab mechanism or continuous shrinkage — before posterior updating, making it computationally and analytically complex. If we could reverse their order, posterior updating would be drastically simplified, for instance, when using a conjugate normal prior on the vector of regression coefficients $\beta$ in the linear regression model \eqref{linear model}. Sparsity can be introduced at the posterior stage via a sparsity-inducing map  
\begin{align} 
\label{sparse projection} 
\beta\mapsto \beta^*:=\arg\min\{\|X\beta-X u\|^2+\lambda_n P(u): u\in \mathbb{R}^p\},
\end{align}
where $P(u)$ is a penalty function such as the $\ell_1$-penalty $\|u\|_1$ and $\lambda_n$ is a tuning parameter. The map is applied to each sample drawn from the (typically conjugate) unrestricted posterior. The induced posterior distribution on $\beta^*$, called the sparse projection-posterior, can be used for inference on $\beta$. Pal and Ghosal \cite{pal2024bayesian} showed that the resulting posterior concentrates near the true regression vector at the optimal rate and is variable selection consistent under the same condition used for the LASSO. Moreover, if a debiasing step is introduced, Bayesian credible intervals for each component have asymptotically correct coverage. They also showed that the sparse projection posterior is a factor of magnitude faster to compute than other Bayesian methods and can be further accelerated using distributed computing. It may be noted that other Bayesian methods for sparse models (except the posterior under a hard-spike-and-slab prior, which is not computable), such as soft-spike-and-slab or continuous shrinkage priors, also implicitly apply a sparsity-inducing step to posterior samples. The main difference is that the step involves a simple thresholding map, whereas \eqref{sparse projection} for the sparse projection-posterior is slightly more involved. However, trading the latter for the former gives the sparse projection posterior the advantage of accessing much faster conjugate posterior draws than MCMC needed for other Bayesian methods. The issue is more fundamental than computational speed, since in the high-dimensional setting, there is no feasible MCMC convergence diagnostic method for the whole chain of $\beta$, since a component of a Markov chain need not be a Markov chain. When $p$ is fixed, Pal and Ghosal \cite{pal2025projection} obtained the weak limit of the sparse projection-posterior random measure that has somewhat like a Bernstein-von Mises type structure. They also showed that a correction on the credibility level obtained from the weak limit achieves asymptotically correct coverage. Pal and Ghosal \cite{pal2025bayesian} extended the approach to address group sparsity and nonparametric additive regression models using a B-spline basis expansion. They derived the optimal posterior contraction and prediction rate under group sparsity, and the consistency of group-variable selection under a version of the beta-min condition. 
The main idea is potentially applicable to many other high-dimensional inference problems, such as generalized regression or graphical models. 

The idea of using the induced posterior through a map that complies with a structural restriction, like sparsity, from a (typically conjugate and easy to sample from) vanilla posterior without the desired structural property has been used in other contexts, such as shape restriction (Chakraborty and Ghosal \cite{chakraborty2021coverage}, Wang and Ghosal \cite{wang2023coverage}) and differential equation models (Bhaumik and Ghosal \cite{bhaumik2015bayesian}). A map that reinforces compliance with the desired structure is often natural and, in many contexts, is given by the projection operator with respect to a metric; in this case, the induced posterior is called the projection posterior. In differential equation models, owing to the similarity with a two-step estimation procedure, the induced posterior was called the two-step posterior in Bhaumik and Ghosal \cite{bhaumik2015bayesian}. In more complex models, such as multivariate monotone function estimation, a projection map may not suffice; an alternative map may be required. Since a structure-compliance-inducing map immerses posterior samples into the desired space, the resulting posterior is called the immersion posterior (Wang and Ghosal \cite{wang2023coverage}), with the projection posterior as a special case. Note that due to the presence of the penalty term in \eqref{sparse projection}, the sparse projection map is not actually a projection map, but is only an immersion map. Pal and Ghosal \cite{pal2024bayesian} gave a Bayesian interpretation of an immersion posterior in terms of a super-parameter model and a misspecified likelihood. 

The immersion posterior approach is similar to VB in that both yield easily computable pseudo-posterior distributions, but not necessarily approximations to a traditional Bayesian posterior. Both are computed via optimization procedures, and theorems that justify their use from a frequentist perspective are available in the literature. On the other hand, the optimization in the VB is performed on distributions, typically via an iterative coordinate descent method, but only needs to be done once. In contrast, for the immersion posterior, the optimization is in the parameter space, which is considerably simpler; the procedure is often non-iterative, but it must be performed for each posterior sample drawn from the unrestricted posterior. The mean-field VB may sometimes underestimate variation, and its credible regions may not have the correct asymptotic coverage, whereas some theorems justifying the immersion posterior uncertainty quantification from the frequentist perspective are available in the literature. 

\subsection{Sparsity priors for nonparametrics and deep learning}

We conclude with a brief section that broadens the discussion to settings in which the object of interest is a function $f$, as is the case, for instance, in nonparametric regression. In such frameworks, using `sparse' methods is also commonplace: for example, thresholding methods estimate the sequence of wavelet coefficients in a basis by setting many empirical coefficients (those not exceeding a certain threshold) to zero. It is therefore particularly natural to consider using spike--and--slab priors $(1-w_n)\delta_0 + w_n \text{Lap}(1)$, for some suitably small $w_n$, in these settings too. In nonparametric white noise regression, Hoffmann et al. \cite{hrs15} proved that putting such priors on the wavelet coefficients of $f$ leads to optimal posterior contraction rates in the supremum norm, rates that are adaptive to the unknown H\"older regularity of the true regression function. In density estimation, spike--and--slab P\'olya trees priors satisfy similar properties, where this time the spike--and--slab $(1-w_n)\delta_{1/2}+w_n\text{Beta}(a,a)$ replaces the usual Beta variables in P\'olya trees  (Castillo and Mismer \cite{cm21}). In this context, one may interpret the spike $\delta_{1/2}$ as a $\text{Beta}(\infty,\infty)$ variable. This idea can also be used to estimate the intensity of point processes, as in Giordano et al. \cite{gkr26}. Nevertheless, in nonconjugate models, as in high-dimensional models, using `hard' spikes can be computationally challenging. A number of recent proposals consider heavy-tailed priors as computationally easier substitutes. Carvalho et al. \cite{ps11} suggested using horseshoe priors, and indeed these can be shown to achieve near-optimal adaptive rates in regression; see Agapiou et al. \cite{ace25}. A simpler alternative is to use heavy-tailed priors with small (but deterministic) scaling factors proposed in Agapiou and Castillo \cite{ac24}, which is both computationally feasible and theoretically (near)-optimal. Here, computations are done via MCMC without a Gibbs step in non-conjugate models. These priors also match theoretical benchmarks for Besov spaces and $\mathbb{L}_r$-risks $(r\ge 2)$; see
 Agapiou et al. \cite{ace25}. 
  
 % The nonparametric regression model where the number of possibly active covariates grows with the number of observations has been studied by Yang and Tokdar \cite{yt15}, who designed a prior based on Gaussian processes where the active coordinates are selected uniformly at random. %It was understood later 
 %Here, a small section on variable selection for GPs and growing dimensions + Bachoc et al. + mentions deep GP discussed below.

In Bayesian deep learning, one main object of interest is again a regression function, this time modeled using a neural network. Sparsity of neural network coefficients plays an important role, both in practice with strategies such as {\em dropout} (Srivastava et al. \cite{dropout}) and in approximation theory for neural networks, where certain sparse networks achieve an optimal approximation-complexity trade-off as demonstrated by Schmidt-Hieber \cite{jsh00}. The next lines extract a (very) small subset of contributions from a rapidly evolving field; we refer readers to Arbel et al. \cite{primer26} for a more focused review on Bayesian deep learning and additional references.

Equipping the network parameters with independent spike-and-slab priors leads to near-optimal posterior rates (see Polson and Ro\v ckov\'a \cite{pr18} and Wang and Ro\v ckov\'a \cite{wr20}). Given the complexity of the output function's dependence on its parameters, an approximation of the posterior is needed for computability. Mean-field Variational Bayes approximations for spike-and-slab priors were considered in Cherief-Abdellatif \cite{badr20}, and Bai et al. \cite{baietal20}.  Continuous shrinkage priors, on the other hand, are particularly natural in this setting. Indeed, in many practical implementations, neural networks are often fully connected (that is, all the weights are non-zero a priori) and overparameterized (that is, the architecture of the network in terms of width and depth is larger than what would in principle be chosen given the approximation--versus--complexity trade-off). In this view, to understand the behavior of such methods, it is natural to consider priors that fit all coefficients, such as continuous shrinkage priors, on overparameterized architectures. Ghosh et al. \cite{vbhs19} considered the implementation of mean-field Variational Bayes using horseshoe priors on coefficients. In Castillo and Egels  \cite{ce25}, a simpler heavy-tailed prior multiplied by a small deterministic scaling parameter was proposed: based on approximation theory for networks, in particular fully connected (see, e.g., Kohler and Langer \cite{kohlerlanger}, Nakada and Imaizumi \cite{imaizumi20}, Suzuki \cite{suzuki21}),  Castillo and Egels \cite{ce25} proved that such a prior on coefficients, based on an overparameterized architecture, leads to a fractional posterior that is fully adaptive to unknown regularity and compositional structures, as well as to data possibly sitting on an unknown small-dimensional manifold. 
%on are natural candidates and achieve, but not necessarily, dense fully connected. 

Yet a different class of priors is that of deep Gaussian processes, introduced in Damianou and Lawrence \cite{dldgp13}, and that are, in their simplest form, compositions of Gaussian processes. Finocchio and Schmidt--Hieber \cite{fsh23} proposed a model-selection-type prior to select active coordinates and showed that it leads to near-minimax posterior rates. Castillo and Randrianarisoa \cite{cr25} showed that the model selection-type prior can be replaced by simple priors on the length-scale parameters of the Gaussian processes, also allowing for a number of possibly active coordinates that grows with the sample size.

\vspace{.5cm}

{{\em Acknowledgments.} We thank Pierre Alquier and Botond Szab\'o, as well as the Associate Editor and three referees, for insightful comments.} 

\section*{Appendix}
\addcontentsline{toc}{section}{Appendix}
\setcounter{section}{9}
\setcounter{subsection}{0}

\subsection{Undirected Graphs}
An undirected graph $G= (V,E)$ consists of a non-empty set of vertices or nodes $V = \{1,\ldots,p\}$ and a set of edges $E \subseteq \{(i,j) \in V \times V: i < j\}$. Nodes connected by an edge are called adjacent nodes. If all the nodes of a graph are adjacent to each other, we have a complete graph. A subgraph $G' = (V',E')$ of $G$, denoted by $G' \subseteq G$ is such that $V' \subseteq V$ and $E' \subseteq E$. A subset $V' \subseteq V$ induces the subgraph $G_{V'} = (V',(V' \times V') \cap E).$ If a subgraph $G_{V'}$ is not contained in any other complete subgraph of $G$, that is, if $G_{V'}$ is a maximal complete subgraph of $G$, then, $V' \subseteq V$ is called a \emph{clique} of $G$. 

A finite subcollection of adjacent edges $(v_0,\ldots,v_{k-1})$ in $G$ forms a path of length $k$. If $v_0 = v_{k-1}$, that is, if the end-points of the path are identical, then we have a $k$-\emph{cycle}. A $\emph{chord}$ of a cycle is a pair of nodes in that cycle which are not consecutive nodes in that path, but are adjacent in $G$. For subgraphs $G_1,G_2$ and $G_3$ of $G$, if every path from a node $v_1 \in G_1$ to a node $v_2 \in G_2$ contains a node in $G_3$, then $G_3$ is said to separate $G_1$ and $G_2$ and is called a \emph{separator} of $G$.

A major portion of our discussion would focus on decomposable graphs. A graph $G$ is said to be {\em decomposable} if every cycle in $G$ having length greater than or equal to four has a chord. The existence of a perfect ordering of the cliques can also characterize the decomposability of a graph. An ordering of the cliques $(C_1,\ldots,C_k) \in \mathcal{C}$ and separators $(S_2,\ldots,S_k) \in \mathcal{S}$ is perfect if it satisfies the \emph{running intersection property}, that is, there exists a $i < j$ such that $S_j = H_{j-1} \cap C_j \subseteq C_i$, where $H_j = \cup_{i=1}^{i}C_j.$ For more details on decomposability and perfect ordering of cliques, we refer the readers to Lauritzen \cite{lauritzen1996graphical}.

\subsection{Graphical models}
An undirected graph $G$ equipped with a probability distribution $P$ on the node-set $V$ is termed an undirected graphical model. We shall focus in particular on Gaussian data, which leads to the concept of Gaussian graphical models. To define the same, we first discuss the Markov property, or the conditional independence property of a $p$-dimensional random variable. A $p$-dimensional random variable ${X} = (X_1,\ldots,X_p)$ is said to be Markov with respect to an undirected graph $G$ if the components $X_i$ and $X_j$ are conditionally independent given the rest of the variables, whenever $(i,j) \notin E$. Thus, an undirected graph $G$ with a Gaussian distribution on the components of the above random variable $X$ corresponding to the nodes $V = \{1,\ldots,p\}$ is called a \emph{Gaussian graphical model} (GGM). Without loss of generality, we can assume that the mean of the Gaussian distribution is zero. For a GGM, there exists a correspondence between the edge set $E$ and the inverse covariance matrix (or, the precision matrix) $\Omega = \Sigma^{-1}$ of $X$ owing to the Markov property. To be precise, whenever $(i,j) \notin E$, the $(i,j)$th element of $\Omega$, given by $\omega_{ij}$, is exactly zero and vice versa. This leads us to consider the cone of positive definite matrices
\begin{equation}
	\mathcal{P}_G = \{\Omega \in \mathcal{M}_p^+: \omega_{ij} = 0, (i,j) \notin E\},
\end{equation}
which defines the parameter space of the precision matrix for a Gaussian distribution defined over a graph $G$. For a decomposable graph $G$, the parameter space for the covariance matrix $\Sigma = (\!(\sigma_{ij})\!)$ is defined by the set $\mathcal{Q}_G$ of partially positive definite matrices, which is a subset of $I_G$, the set of incomplete matrices with missing entries $\sigma_{ij}$ whenever $(i,j) \notin E$. Then,
\begin{equation}
	\mathcal{Q}_G = \{B \in \mathcal{I}_G: B_{C_i} > 0,\; i = 1,\ldots,k\}. 
\end{equation}
Gr{\"o}ne et al. \cite{grone1984positive} showed that there is a bijection between the spaces $\mathcal{P}_G$ and $\mathcal{Q}_G$. To be precise, for $\Omega \in \mathcal{P}_G$, we can define $\mathcal{Q}_G$ as the parameter space for the GGM where $\Sigma = \kappa(\Omega^{-1})$, and $\kappa:\mathcal{M}_p^+ \rightarrow \mathcal{I}_G$ is the projection of $\mathcal{M}_p^+$ into $I_G$. Thus, a GGM Markov with respect to a graph $G$ is given by the family of distributions
\begin{equation}
	\mathcal{N}_G = \{\mathrm{N}_p({0},\Sigma), \Sigma \in \mathcal{Q}_G\} = \{\mathrm{N}_p({0},\Omega^{-1}), \Omega \in \mathcal{P}_G\}.
\end{equation}

\subsection{Hyper Inverse-Wishart and G-Wishart priors}

The inverse-Wishart distribution $\mathrm{IW}_p(\delta,D)$ with degrees of freedom $\delta$  and a $p$-dimensional positive definite fixed scale matrix $D$ is the conjugate prior for the covariance matrix $\Sigma$ in the case of a complete graph. We denote $\Sigma \sim \mathrm{IW}_p(\delta,D)$ having density
\begin{equation}
	p(\Sigma \mid \delta, D) = c(\delta,p)\{\det(\Sigma)\}^{-(\delta + 2p)/2}\exp\{-\mathrm{tr}(\Sigma^{-1}D)/2\},
\end{equation}
where $c(\delta,p) = \{\det(D)/2\}^{(\delta + p - 1)/2}/{\Gamma_p\{(\delta+p-1)/2\}}$ is the normalizing constant and  $\Gamma_p(t) = \pi^{p(p-1)/4}\prod_{j=1}^{p}\Gamma(t - (j-1)/2)$ is the $p$-variate gamma function. 

For decomposable graphs, recall that a perfect set of cliques always exists.  In that case, we can write the density as
\begin{equation}
	p(x \mid \Sigma, G) \propto \frac{\prod_{C \in \mathcal{C}}p_C({x}_C \mid \Sigma_C)}{\prod_{S \in \mathcal{S}}p_S({x}_S \mid \Sigma_S)},
\end{equation}
which is the Markov ratio of respective marginal Gaussian distributions corresponding to the cliques and separators of $G$. Dawid and Lauritzen \cite{dawid1993hyper} came up with a generalization of the inverse-Wishart, called the hyper inverse-Wishart, which is conjugate to the Gaussian distribution Markov with respect to the graph $G$. The form of the prior distribution depends on the clique-marginal covariance matrices $\Sigma_C$, such that $\Sigma_C \sim \mathrm{IW}_{|C|}(\delta,D_C)$, where $D_C$ is the submatrix of the scale matrix $C$ induced by the clique $C \in \mathcal{C}.$ The hyper inverse-Wishart is thus constructed on the parameter space $\mathcal{Q}_G$ with density given by
\begin{equation}
	p_G(\Sigma \mid \delta, D) = \frac{\prod_{C \in \mathcal{C}}p(\Sigma_C \mid \delta, D_C)}{\prod_{S \in \mathcal{S}}p(\Sigma_S \mid \delta, D_S)},
\end{equation}
where $p(\cdot)$ refers to the density of the respective inverse-Wishart distribution.

For a complete graph, the inverse-Wishart prior on the covariance matrix induces a conjugate prior for the precision matrix $\Omega$, namely, the Wishart distribution, with density
\begin{equation}
	p(\Omega \mid \delta, D) = c(\delta,p)\{\det(\Omega)\}^{(\delta - 2)/2}\exp\{-\mathrm{tr}(\Omega D)/2\},
\end{equation} 
with an identical normalizing constant as in inverse-Wishart. Similar in lines with this, for a decomposable graphical model, the hyper inverse-Wishart prior induces a prior distribution on the precision matrix $\Omega$ with density
\begin{equation}
	p_G(\Omega \mid \delta,D) = I_G(\delta,D)^{-1}\{\det(\Omega)\}^{(\delta-2)/2}\exp\{-\mathrm{tr}(D\Omega)/2\},
\end{equation}
where $I_G(\delta,A)$ is the normalizing constant given by
\begin{equation}
	\label{G-Wishart normalization}
	I_G(\delta,D) = \frac{\prod_{S \in \mathcal{S}}\{\det(D_S)\}^{(\delta+|S|-1)/2}/\Gamma_{|S|}({(\delta+|S|-1)}/{2})}{\prod_{C \in \mathcal{C}}\{\det(D_C)\}^{(\delta+|C|-1)/2}/\Gamma_{|C|}({(\delta+|C|-1)}/{2})}.
\end{equation}
This distribution on $\Omega$ is called the G-Wishart distribution $\mathrm{W}_G$. It forms the Diaconis-Ylvisaker conjugate prior for the precision matrix of a Gaussian distribution Markov with respect to the decomposable graph $G$. However, unlike the hyper inverse-Wishart prior, the G-Wishart prior can be generalized to non-decomposable graphical models as well, although there is no closed-form analytical expression of the normalizing constant $I_G(\delta,D)$, except for a result by Uhler et al. \cite{uhler2018exact}. 

Letac and Massam \cite{letac2007wishart} introduced a more general class of conjugate priors, namely the $\mathrm{W}_{\mathcal{P}_G}$-Wishart class of priors for the precision matrix $\Omega$ with three sets of parameters ---  ${\alpha}=(\alpha_1,\ldots,\alpha_p)$ and ${\beta}=(\beta_1,\ldots,\beta_p)$ which are suitable functions defined on the cliques and separators of the graph, and a scale matrix $D$. We note that the above class of prior distributions is the Wishart distribution for a fully connected graph and includes the G-Wishart distribution for decomposable graphs as a special case with suitable choices of ${\alpha}$ and ${\beta} $.

\subsection{Posterior computation in graphical models}

\subsubsection{Exact expressions under conjugacy}

Denoting by $\mathbb{X}$ the $n \times p$ data matrix corresponding to a random sample from a Gaussian distribution $\mathrm{N}_p({0},\Omega^{-1})$, conditional on a specific graph $G$, the posterior density $p(\Omega \mid \mathbb{X}, G)$  of $\Omega$ is given by
\begin{equation}
	I_G(\delta+n, D + S_n)^{-1}\{\det(\Omega)\}^{(\delta+n-2)/2}\exp[-\mathrm{tr}\{(D+S_n)\Omega\}/2 ],
\end{equation}
where $S_n = \mathbb{X}\mathbb{X}^T.$

Carvalho et al. \cite{carvalho2007simulation} proposed a sampler for the hyper inverse-Wishart distribution corresponding to a decomposable graph based on distributions of submatrices of the covariance matrix $\Sigma$. Samples drawn from the hyper-inverse-Wishart distribution can be inverted to obtain samples from the corresponding G-Wishart distribution.

Rajaratnam et al. \cite{rajaratnam2008flexible} provided a closed form expression of the posterior mean $\mathrm{E}(\Omega \mid S_n)$ of $\Omega$ with a $\mathrm{W}_{\mathcal{P}_G}$-Wishart prior as
\begin{equation}
	-2\Big[ \sum_{j=1}^{k}\big(\alpha_j - \frac{n}{2}\big)\big(\big(D + \kappa(nS_n)\big)_{C_j}^{-1}\big)^0  - \sum_{j=2}^{k}\big(\beta_j - \frac{n}{2}\big)\big(\big(D + \kappa(nS_n)\big)_{S_j}^{-1}\big)^0\Big],
\end{equation}
where for a $p\times p$ matrix $A = (\!(a_{ij})\!)$, $A_T^{-1}$ denotes the inverse $(A_T)^{-1}$ of the submatrix $A_T$ and $(A_T)^0 = (\!(a^*_{ij})\!)$ denotes a $p$-dimensional matrix such that $a^*_{ij} = a_{ij}$ for $(i,j) \in T \times T$, and $0$ otherwise. Here $\alpha_j$ and $\beta_j$ are the $j$th component of the parameters ${\alpha}$ and ${\beta}$ respectively, $j=1,\ldots,p$. They also obtained the expression for the Bayes estimator with respect to Stein's loss 
$L(\hat{\Omega},\Omega) = \mathrm{tr}(\hat{\Omega} - \Omega)^2$ under the $\mathrm{W}_{\mathcal{P}_G}$-Wishart prior.

\subsubsection{MCMC sampling methods for G-Wishart prior}

Madigan and York \cite{madigan1995bayesian} proposed an approach based on the Metropolis-Hastings algorithm to traverse the space of decomposable graphs $\mathcal{G}$, while Giudicci and Green \cite{giudici1999decomposable} used a reversible-jump MCMC sampler that also included the precision matrix as one of the state variables.

When the graph is not necessarily decomposable, the G-Wishart prior is still conjugate, but the normalizing constant $I_G(\delta,D)$ does not have a simple expression. Only recently, Uhler et al. \cite{uhler2018exact} derived an expression, but computing it for inference remains challenging. Lenkoski and Dobra \cite{lenkoski2011computational} developed a Laplace approximation 
$$
\hat{I_G(\delta,D)} = \exp\{l(\hat{\Omega})\}(2\pi)^{p/2}\det(H(\hat{\Omega}))^{-1/2},
$$
where $l = ((\delta - 2)\log \det(\Omega) - \mathrm{tr}(D\Omega))/2$, the matrix $\hat{\Omega}$ is the mode of the G-Wishart density and $H$ is the Hessian matrix associated with $l$, but it lacks accuracy unless the clique sizes are small. 

Atay-Kayis and Massam \cite{atay2005monte} proposed a Monte Carlo-based method to sample from the G-Wishart distribution as well as to compute the prior normalizing constant. They considered the Cholesky decomposition $D^{-1} = Q^{\mathrm{T}} Q$ and $\Omega = \Phi^{\mathrm{T}}\Phi$, and defined $\Psi = \Phi Q$. Since for $(i,j) \notin E$, $\omega_{ij}=0$, and hence the elements $\psi_{ij}$ of $\Psi$ for $(i,j) \notin E$ are functions of $\psi_{ij}$, $(i,j) \in E$ and $\psi_{ii}$, $i=1,\ldots,p$. Thus the free elements appearing in the Cholesky decomposition of $\Omega$ are given by $\Psi^E = \left(\psi_{ij},\,(i,j) \in E;\, \psi_{ii},i=1,\ldots,p\right)$. They showed that the free elements have a density
$$
p(\Psi^E) \propto f(\Psi^E)h(\Psi^E),
\qquad f(\Psi^E) = \exp\{-\sum_{(i,j) \notin E}\psi_{ij}^2/2\},$$ 
is a function of the non-free elements of $\Psi$, which in turn can be uniquely expressed in terms of $\Psi^E$. Furthermore, $h(\Psi^E)$ is the product of densities of random variables $\psi_{ii}^2 \sim \chi_{\delta + v_i}^2,\,i=1,\ldots,p,$ and $\psi_{ij} \sim \mathrm{N}(0,1),\, (i,j) \in E$, where $v_i$ stands for the cardinality of $\{j: j > i, (i,j) \in E\}.$ Generating samples for $\Omega$ then uses an acceptance-rejection method based on the above density function of the free elements $\Psi^E$. The normalizing constant $I_G(\delta,D)$ can be expressed as the product of a known constant and the expected value of $f(\Psi^E)$. Using samples from the distribution of $\Psi^E$, straightforward Monte Carlo computation gives $I_G(\delta,D)$. However, the Monte Carlo integration method is computationally expensive for non-complete large prime components of the graph, owing to a matrix completion step involving the non-free elements of $\Psi$. 

To alleviate this problem, Wang and Carvalho \cite{wang2010simulation} used a prime-component decomposition so that sampling is individually carried out in lower-dimensional subgraphs of $G$. However, owing to the dependencies of the sampler on $(\delta,D,G)$, the acceptance rate of the MCMC procedure can be very low. Mitsakakis et al. \cite{mitsakakis2011metropolis} developed an independent Metropolis-Hastings algorithm for sampling from $\mathrm{W}_G(\delta,D)$ based on the density of $\Psi^E$. This method, though, improves the acceptance rate over that proposed by Wang and Carvalho \cite{wang2010simulation}, suffers from low acceptance rates and slow mixing in large graphs. In line with Mitsakakis et al. \cite{mitsakakis2011metropolis}, Dobra et al. \cite{dobra2011bayesian} used a random walk Metropolis-Hastings sampler, in which only one entry of $\Psi^E$ is perturbed in a single step, rather than changing all entries of $\Psi^E$ as in the former approaches. Although this method improves the sampler's efficiency, it still suffers from the matrix completion problem for non-free elements, incurring a time complexity of $O(p^4)$ per Monte Carlo sample. This results in painfully slow computations for large graphs.

Jones et al. \cite{jones2005experiments} used the method of Atay-Kayis and Massam \cite{atay2005monte} to compute prior and posterior normalizing constants of the G-Wishart distribution to traverse the space of graphs using a technique called the \emph{Stochastic Shotgun Search Algorithm} using the uniform or Erd\"os-Renyi type priors on the space of graphs. Though the method performs well in low dimensions, it fails to scale up in high dimensions due to the huge search space.

Lenkoski and Dobra \cite{lenkoski2011computational} and Mitsakakis et al. \cite{mitsakakis2011metropolis} proposed to use the Bayesian iterative proportional scaling algorithm developed by Piccioni \cite{piccioni2000independence} and Asci and Piccioni \cite{asci2007functionally} to sample from the G-Wishart distribution. Their method requires enumerating the maximum cliques of $G$, an NP-hard problem, and inverting large matrices, both of which are computationally burdensome.

As discussed before, Giudicci and Green \cite{giudici1999decomposable} developed a reversible-jump MCMC method to learn the graphical structure by sampling over the joint space for $(\Omega,G)$. Dobra et al. \cite{dobra2011bayesian} also developed a reversible jump MCMC method over $(\Omega,G)$, thus avoiding the issue of searching over a large space of graphs. But their method cannot avoid problems with the computation of prior normalizing constants and matrix completion. The crucial bottleneck in computing the acceptance probabilities of the MCMC based procedures developed for joint exploration of $(\Omega,G)$ is the ratio of prior normalizing constants $I_{G^{-e}}(\delta,D)/I_G(\delta,D),$ where $G^{-e}$ is the graph obtained from $G = (V,E)$ by deleting a single edge $e \in E.$ For the choice of $D = I_p$, Letac et al. \cite{letac2017ratio} showed that the above ratio can be reasonably approximated by constant times the ratio of two gamma functions as
\begin{equation}
	{I_{G^{-e}}(\delta,D)}/{I_G(\delta,D)} \approx ({2\sqrt{\pi}})^{-1}{\Gamma\left(\frac{\delta + d}{2}\right)}/{\Gamma\left(\frac{\delta+d+1}{2}\right)},
\end{equation}
where $d$ is the number of paths of length two between the two nodes in edge $e$. They showed that under certain conditions and graph configurations, the approximation performs reasonably well.

To circumvent the problems of computing the normalizing constant $I_G(\delta,D)$ for arbitrary graphs, especially in the high-dimensional scenario, Wang and Li \cite{wang2012efficient} proposed a double reversible jump MCMC algorithm using the partial analytic structure (PAS) of the G-Wishart distribution. Their proposed algorithm involves a Metropolis-Hastings step to move from $G = (V,E)$ to $G' = (V,E')$, where $E'$ differs in only one edge from $E$, with the acceptance probability for the Metropolis-Hastings step having exact analytical expressions which are easy to compute. 

A completely new approach, called the \emph{Birth-and-Death MCMC} method (BDMCMC) for graphical model selection, was introduced by Mohammadi and Wit \cite{mohammadi2015bayesian} and explores the space of graphs by adding (birth) or deleting (death) an edge. They determine a continuous time birth-death Markov process on $\Omega$, more specifically, an independent Poisson process on $\Omega$ such that under suitable conditions, the process has a stationary distribution $p(\Omega,G \mid {X}^{(n)})$. The birth and death rates are given by
\begin{eqnarray}
	\beta_e(\Omega) &=& \frac{P(G^{+e},\Omega^{+e}\backslash (\omega_{ij},\omega_{jj}) \mid {X}^{(n)})}{P(G,\Omega \backslash \omega_{jj} \mid {X}^{(n)})},\quad \mbox{ for each } e \in \bar{E}, \nonumber \\
	\delta_e(\Omega) &=& \frac{P(G^{-e},\Omega^{-e}\backslash \omega_{jj} \mid {X}^{(n)})}{P(G,\Omega \backslash (\omega_{ij},\omega_{jj}) \mid {X}^{(n)})},\quad \mbox{ for each } e \in E; 
\end{eqnarray}
here $G^{+e},G^{-e}$ are graphs obtained from $G$ after addition or removal of an edge $e$ from $G$. The direct-sampler method proposed by Lenkoski \cite{lenkoski2013direct} is used to sample from the posterior distribution of $\Omega$. The BDMCMC algorithm is easy to implement, scales well to high-dimensional problems, and performs better at structure learning than the Bayesian methods discussed earlier. 

In a recent work, Bhadra et al. \cite{bhadra2024evidence} proposed a novel approach based on a telescoping block decomposition of the precision matrix for calculating marginal likelihood (also known as \emph{model evidence}) of graphical models under various priors, including the $G$-Wishart and graphical shrinkage priors. In addition to addressing the open problem of computing model evidence for graphical models, this method also led to designing a column-wise sampler for the $G$-Wishart, which is significantly faster than the direct sampler method of Uhler et al. \cite{lenkoski2013direct}.

\subsubsection{MCMC sampling methods for the Bayesian graphical LASSO and variants}

Wang \cite{wang2012bayesian} developed a block Gibbs sampler to simulate from the posterior distribution of the Bayesian graphical LASSO through a reparameterization $T = (\!(\tau_{ij})\!)$ with zeros as diagonal entries and $\tau$ in upper diagonal entries. Then, the matrices $\Omega, S$, and $T$ are partitioned with respect to the last column and row as
\begin{equation}
	\Omega = \left(\begin{array}{cc}
		\Omega_{11}& {\omega}_{12}  \\ 
		{\omega}_{12}^T& \omega_{22} 
	\end{array} \right) ,\;
	S = \left(\begin{array}{cc}
		S_{11}& {s}_{12}  \\ 
		{s}_{12}^T& s_{22} 
	\end{array} \right), \; 
	\Upsilon = \left(\begin{array}{cc}
		\Upsilon_{11}& {\tau}_{12}  \\ 
		{\tau}_{12}^T& 0 
	\end{array} \right) .
\end{equation}
With the change of variable $(\omega_{12},\omega_{22}) \mapsto  (\beta = \omega_{12},\gamma = \omega_{22} - \omega_{12}^T\Omega_{11}^{-1}\omega_{12}),$
the conditional posterior distribution of $(\beta, \gamma)$ given the rest is given by
\begin{eqnarray*}
	p(\beta,\gamma, \Omega_{11},T,{X}^{(n)},\lambda) &\propto & \gamma^{n/2}\exp\Big(\frac{s_{22} + \lambda}{2}\gamma \Big) \nonumber \\
	&& \times \exp\Big(-\frac{1}{2}[\beta^T\{D_{\tau}^{-1} + (s_{22}+ \lambda)\Omega_{11}^{-1}\}\beta + 2s_{12}^T\beta] \Big),
\end{eqnarray*}
where $D_{\tau} = \mathrm{diag}(\tau_{12}),$ so that the conditional distributions of $\gamma$ and $\beta$ are independent gamma and a normal, respectively, whereas the inverse of the latent scale parameters $\tau_{ij}$ are independent inverse-Gaussian distributions. The resulting block Gibbs sampler iteratively samples one column of $\Omega$ at a time. It is interesting to note that the positive definiteness constraint is also maintained throughout because $\gamma$ is always positive definite. The posterior samples hence obtained cannot be directly used for structure learning, though, since the Bayesian graphical LASSO prior puts zero mass on the event $\{\omega_{ij} = 0\}.$  

As discussed in the previous sections, computation poses significant challenges in Bayesian graphical models for high-dimensional settings. To deal with large dimensions in arbitrary graphs, Wang \cite{wang2015scaling} developed a new technique called \emph{stochastic search structure learning} for precision matrices as well as covariance matrices in a graphical model setup, based on soft-spike-and-slab priors on the elements of the matrix. This work also focuses on the issue of structure learning of graphs using a fully Bayesian model by specifying priors on binary variables ${Z} = (z_{ij})_{i<j}$, which are indicators for the edges of the underlying graph $G$. The hierarchical prior is then specified as
\begin{eqnarray}
	p(\Omega \mid {Z}, \theta) &=& C({Z},\nu_0,\nu_1,\lambda)^{-1}\prod_{i < j}\mathrm{N}(\omega_{ij} \mid 0, \nu_{z_{ij}}^2)\prod_{i=1}^{p}\mathrm{Exp}(\omega_{ii} \mid \lambda/2), \nonumber \\
	p({Z}\mid \theta) &=& C(\theta)^{-1}C({Z},\nu_0,\nu_1,\lambda)\prod_{i < j}\pi^{z_{ij}}(1-\pi)^{z_{ij}},
\end{eqnarray}
where  $\theta = (\nu_0,\nu_1,\pi,\lambda),$ $\mathrm{N}(a\mid 0,\nu^2)$ is the density of a normal distribution with mean zero and variance $\nu^2$, $C(\theta)$ and $C({Z},\nu_0,\nu_1,\lambda)$ are normalizing constants, and $\nu_{z_{ij}}$ is $\nu_0$ or $\nu_1$ according as $z_{ij}$ is 0 or 1. The hyperparameter $\pi \in (0,1)$ controls the prior distribution of the binary edge indicators in ${Z}$. The hierarchical prior specification above leads to the following prior specification on $\Omega$:
\begin{eqnarray}
	p(\Omega) &=& C(\theta)^{-1}\prod_{i < j}\left\{(1-\pi)\mathrm{N}(\omega_{ij} \mid 0, \nu_0^2) + \pi \mathrm{N}(\omega_{ij} \mid 0, \nu_1^2) \right\} \nonumber \\
	&& \qquad\times \prod_{i=1}^{p}\mathrm{Exp}(\omega_{ii} \mid \lambda/2) \mathbbm{1}(\Omega \in \mathcal{M}_p^+),
\end{eqnarray} 
A small $\nu_0 > 0$ and a large value of $\nu_1 > 0$ induce a soft-spike-and-slab prior on the elements of $\Omega$. The above prior, having a two-component mixture of normals, facilitates graph structure learning via the latent binary indicators ${Z}$. The sampling procedure for generating samples from the posterior distribution $p(\Omega \mid {X}^{(n)},{Z})$ is identical to the block Gibbs sampler proposed in Wang \cite{wang2012bayesian}, via introduction of the $p$-dimensional symmetric matrix ${V} = (\!(\nu_{z_{ij}}^2)\!)$ with zeros in diagonal and $(\!(\nu_{ij}^2: {i<j})\!)$ as the upper-diagonal entries, where $\nu_{ij}^2 = \nu_{z_{ij}}^2$. The conditional posterior distributions of the binary indicator variables ${Z}$ given the data and $\Omega$ are independent Bernoulli with success probability
\begin{equation}
	p(z_{ij} = 1 \mid \Omega, {X}^{(n)}) = \frac{\pi \mathrm{N}(\omega_{ij} \mid 0, \nu_1^2)}{\pi \mathrm{N}(\omega_{ij} \mid 0, \nu_1^2) + (1-\pi) \mathrm{N}(\omega_{ij} \mid 0, \nu_0^2)}.
\end{equation}
Although the structure-learning accuracy of the above method is comparable to that of other Bayesian structure-learning methods, the block Gibbs sampler-based computational approach makes it simpler and faster, especially when scaling to large dimensions. 

\subsubsection{Laplace approximation to compute posterior model probabilities}

For graphical structure learning in the absence of explicit expressions, one may compute the posterior probabilities of various models using reversible-jump Markov chain Monte Carlo, which is computationally expensive. Banerjee and Ghosal \cite{banerjee2015bayesian} proposed directly computing the marginal posterior probabilities of models using a Laplace approximation. The main idea is to expand the log-likelihood around the posterior mode, which can be easily identified as the graphical LASSO within that model and computed efficiently, and then integrate the resulting approximate likelihood function. The resulting computation is speedy, but a drawback is that the approach works only for `regular models', where none of the components of the posterior mode is $0$, because the log-likelihood function is singular at any point having a coordinate $0$. The problem is partly alleviated by the fact that for every `non-regular model', there is a regular submodel with a higher posterior probability. Hence, at least for model selection, the approach can be restricted to the latter class.  

\bibliographystyle{imsart-number}
\bibliography{bibreview_24}
\end{document}